\documentclass[11pt,centertags,reqno,twoside]{amsart}
\usepackage{amssymb}
\usepackage{mathptmx}

\DeclareSymbolFont{SY}{U}{psy}{m}{n}
\DeclareMathSymbol{\emptyset}{\mathord}{SY}{'306}

\renewcommand{\eqref}[1]{{\rm(\ref{#1})}}

\newcommand{\bbC}{{\mathbb C}}
\newcommand{\bbR}{{\mathbb R}}

\newcommand{\cB}{{\mathcal B}}

\newcommand{\cG}{{\mathcal G}}

\newcommand{\conv}{\mathop{\rm conv}}
\newcommand{\cO}{{\mathcal O}}
\newcommand{\cP}{{\mathcal P}}

\newcommand{\cT}{{\mathcal T}}

\newcommand{\ri}{{\rm i}}
\newcommand{\sA}{{\sf A}}
\newcommand{\sB}{{\sf B}}

\newcommand{\sE}{{\sf E}}

\newcommand{\sL}{{\sf L}}

\newcommand{\sS}{\text{\sffamily{\textit{S}}}}
\newcommand{\sV}{{\sf V}}

\newcommand{\fG}{\mathfrak{G}}
\newcommand{\fH}{\mathfrak{H}}
\newcommand{\fK}{\mathfrak{K}}
\newcommand{\fL}{\mathfrak{L}}
\newcommand{\fM}{\mathfrak{M}}
\newcommand{\fN}{\mathfrak{N}}

\newcommand{\diag}{\mathop{\rm diag}}
\newcommand{\dist}{\mathop{\rm dist}}
\newcommand{\Img}{\mathop{\rm Im}}

\newcommand{\lal}{{\langle}}
\newcommand{\ral}{{\rangle}}

\newcommand{\be}{\begin{equation}}
\newcommand{\ee}{\end{equation}}

\DeclareMathOperator{\spec}{spec}

\newcommand{\ran}{\mathop{\mathrm{Ran}}}
\newcommand{\Ran}{\mathop{\mathrm{Ran}}}
\newcommand{\dom}{\mathop{\mathrm{Dom}}}
\newcommand{\Dom}{\mathop{\mathrm{Dom}}}

\allowdisplaybreaks

\numberwithin{equation}{section}

\newtheorem{theorem}{Theorem}[section]
\newtheorem{corollary}[theorem]{Corollary}
\newtheorem{lemma}[theorem]{Lemma}

\newtheorem{hypothesis}[theorem]{Hypothesis}
\theoremstyle{definition}
\newtheorem{definition}[theorem]{Definition}
\newtheorem{remark}[theorem]{Remark}
\newtheorem{example}[theorem]{Example}

\begin{document}

\title[Bounds on variation of spectral subspaces under
$J$-self-adjoint perturbations]
{Bounds on variation of spectral subspaces \\ under
$J$-self-adjoint perturbations$^*$}\thanks{$^*$This work was supported
by the Deutsche For\-sch\-ungs\-gemeinschaft (DFG), the
Heisenberg-Landau Program, and the Russian Foundation for Basic
Research.}

\author[S. Albeverio, A. K. Motovilov, and A. A. Shkalikov]
{Sergio Albeverio,  Alexander K. Motovilov, and Andrei A. Shkalikov}

\address{Sergio Albeverio,
Institut f\"ur Angewandte Mathematik, Universit\"at Bonn,
Wegelerstra{\ss}e 6, D-53115 Bonn, Germany; SFB 611 and HCM, Bonn; BiBoS,
Bielefeld-Bonn; CERFIM, Lo\-car\-no; Accademia di Architettura, USI, Mendrisio}
\email{albeverio@uni-bonn.de}

\address{Alexander K. Motovilov, Bogoliubov Laboratory of
Theoretical Physics, JINR, Joliot-Cu\-rie 6, 141980 Dubna, Moscow
Region, Russia} \email{motovilv@theor.jinr.ru}

\address{Andrei A. Shkalikov, Faculty of Mathematics and Mechanics,
Moscow Lomonosov State University, Le\-nin\-skie Gory, 119992 Moscow, Russia}
\email{ashkalikov@yahoo.com}




\begin{abstract}
Let $A$ be a self-adjoint operator on a Hilbert space $\fH$. Assume
that the spectrum of $A$ consists of two disjoint components
$\sigma_0$ and $\sigma_1$. Let $V$ be a bounded operator on $\fH$,
off-diagonal and $J$-self-adjoint with respect to the orthogonal
decomposition $\fH=\fH_0\oplus\fH_1$ where $\fH_0$ and $\fH_1$ are
the spectral subspaces of $A$ associated with the spectral sets
$\sigma_0$ and $\sigma_1$, respectively. We find (optimal)
conditions on $V$ guaranteeing that the perturbed operator $L=A+V$
is similar to a self-adjoint operator. Moreover, we prove a number
of (sharp) norm bounds on the variation of the spectral subspaces of $A$
under the perturbation $V$. Some of the results obtained are
reformulated in terms of the Krein space theory. As an example, the
quantum harmonic oscillator under a $\cP\cT$-symmetric perturbation
is discussed.
\end{abstract}

\maketitle

\section{Introduction}
\label{SIntro}
Let $A$ be a (possibly unbounded) self-adjoint operator on a Hilbert
space $\fH$. Assume that $V$ is a bounded operator on $\fH$. It
is well known that in such a case the spectrum of the perturbed
operator $L=A+V$ lies in the closed $\|V\|$-neighborhood of the
spectrum of $A$ even if $V$ is non-self-adjoint. Thus, if the
spectrum of $A$ consists of two disjoint components $\sigma_0$ and
$\sigma_1$, that is,~if
\begin{equation}
\label{ddist}
\spec(A)=\sigma_0\cup\sigma_1\text{\, and \,}
\dist(\sigma_{0}, \sigma_{1})=d  > 0,
\end{equation}
then the perturbation $V$ with a sufficiently small norm does not
close the gaps between $\sigma_0$ and $\sigma_1$ in $\bbC$. This
allows one to think of the corresponding disjoint spectral
components $\sigma'_{0}$ and $\sigma'_{1}$ of the perturbed operator
$L=A + V$ as a result of the perturbation of the spectral sets
$\sigma_{0}$ and $\sigma_{1}$, respectively.

Assuming \eqref{ddist}, by $\sE_A(\sigma_{0})$ and
$\sE_A(\sigma_{1})$ we denote the spectral projections of $A$
associated with the disjoint Borel sets $\sigma_0$ and $\sigma_1$,
and by $\fH_0$ and $\fH_1$ the respective spectral subspaces,
$\fH_0=\ran \sE_A(\sigma_{0})$ and $\fH_1=\ran\sE_A(\sigma_{1})$. If
there is a possibility to associate with the disjoint spectral sets
$\sigma'_0$ and $\sigma'_1$ the corresponding spectral subspaces of
the perturbed (non-self-adjoint) operator $L=A+V$, we denote them by
$\fH'_0$ and $\fH'_1$. In particular, if one of the sets $\sigma'_0$
and $\sigma'_1$ is bounded, this can easily be done by using the Riesz
projections (see, e.g. \cite[Sec.~III.4]{Kato}).

In the present note we are mainly concerned with bounded
perturbations $V$ that possess the property
\begin{equation}
\label{VJJ}
V^*=JVJ,
\end{equation}
where $J$ is a self-adjoint involution on $\fH$ given by
\begin{align}
\label{Jinv}
J=\sE_A(\sigma_{0})-\sE_A(\sigma_{1}).
\end{align}
Operators $V$ with the property \eqref{VJJ} are called $J$-self-adjoint.

A bounded perturbation $V$ is called diagonal with respect to the
orthogonal decomposition
$\fH=\fH_0\oplus\fH_1$
if it commutes with the involution $J$, $V\!J=JV$. If $V$
anticommutes with $J$, i.e. $V\!J=-JV$, then $V$ is said to be
off-diagonal. Clearly, any bounded $V$ can be represented as the sum
$V=V_\mathrm{diag}+V_\mathrm{off}$ of the diagonal,
$V_\mathrm{diag}$, and off-diagonal, $V_\mathrm{off}$, terms. The
spectral subspaces $\fH_0$ and $\fH_1$ remain invariant under
$A+V_\mathrm{diag}$ while adding a non-zero $V_\mathrm{off}$ does
break the invariance of $\fH_0$ and $\fH_1$. Thus, the core of the
perturbation theory for spectral subspaces is in the study of their
variation under off-diagonal perturbations (cf. \cite{KMM2}). This
is the reason why we add to the hypothesis \eqref{VJJ} another basic
assumption, namely that all the perturbations $V$ involved are
off-diagonal with respect to the decomposition
$\fH=\fH_0\oplus\fH_1$.

We recall that if an off-diagonal perturbation $V$ is self-adjoint
in the usual sense, that is, $V^*=V$, then the condition
\begin{equation}
\label{Vd2I}
\|V\|<\frac{d}{2}
\end{equation}
ensuring the existence of gaps between the perturbed spectral sets
$\sigma'_0$ and $\sigma'_1$ may be essentially relaxed. Generically,
if no assumptions on the mutual position of the initial spectral
sets $\sigma_0$ and $\sigma_1$ are made except \eqref{ddist}, the
sets $\sigma'_0$ and $\sigma'_1$ remain disjoint for any
off-diagonal self-adjoint $V$ satisfying the bound
$\|V\|<\frac{\sqrt{3}}{2}d$ (see \cite[Theorem 1 (i)]{KMM4}). If, in
addition to \eqref{ddist}, it is known that one of the sets
$\sigma_0$ and $\sigma_1$ lies in a finite gap of the other set then
this bound may be relaxed further: for the perturbed sets
$\sigma'_0$ and $\sigma'_1$ to be disjoint it only suffices to
require that $\|V\|<\sqrt{2}d$ (see \cite[Theorem 2 (i)]{KMM4}; cf.
\cite[Remark 3.3]{KMM3}). Finally, if the sets $\sigma_0$ and
$\sigma_1$ are subordinated, say $\sup\sigma_0<\inf\sigma_1$, then
no requirements on $\|V\|$ are needed at all: the interval
$(\sup\sigma_0,\inf\sigma_1)$ belongs to the resolvent set of the
perturbed operator $L=A+V$ for any bounded off-diagonal self-adjoint
$V$ (see \cite{AdLT,DK70,LT2001}; cf. \cite{KMM5}) and even for some
off-diagonal unbounded symmetric $V$ (see \cite[Theorem 1]{MotSel}).
It is easily seen from Example \ref{ExNS} below that in the case of
$J$-self-adjoint off-diagonal perturbations the condition
\eqref{Vd2I} ensuring the disjointness of the perturbed spectral
sets $\sigma'_0$ and $\sigma'_1$ can be relaxed for none of the
above dispositions of the initial spectral sets $\sigma_0$ and
$\sigma_1$.

Assuming that $V$ is a bounded $J$-self-adjoint off-diagonal
perturbation of the (possibly un\-bo\-un\-ded) self-adjoint ope\-ra\-tor
$A$ we address the following questions:
\begin{enumerate}
\item[(i)] Does the spectrum of the perturbed operator $L=A+V$ remain real under
conditions \eqref{ddist} and \eqref{Vd2I}?

\item[(ii)] If yes, is it then true that $L$ is similar to a self-adjoint operator?
\item[(iii)] What are the (sharp) bounds on variation of the spectral subspaces associated
with the spectral sets $\sigma_0$ and $\sigma_1$ as well as on the
variation of these sets themselves?
\end{enumerate}

In our answers to the above questions we distinguish two cases:
\begin{enumerate}
\item[(G)] the generic case where
no assumptions on the mutual positions of the spectral sets
$\sigma_0$ and $\sigma_1$ are done except for the disjointness
assumption \eqref{ddist};

\item[(S)] the particular case where the sets $\sigma_0$ and
$\sigma_1$ are either subordinated, e.g.
$\sup\sigma_0<\inf\sigma_1$, or one of these sets lies in a finite
gap of the other set, say $\sigma_0$ lies in a finite gap of
$\sigma_1$.
\end{enumerate}

We have to underline that this distinction is quite different from
the one that arises when the perturbations $V$ are self-adjoint in
the usual sense: the case (S) now combines the two spectral
dispositions that should be treated separately if $V$ were
self-adjoint (see \cite{AlMoSel,DK70,KMM4,MotSel}).

Our answers to the questions (i) and (ii) are complete and positive
in the case (S). In this case the spectrum of the perturbed operator
$L=A+V$ does remain real for any off-diagonal $J$-self-adjoint $V$
satisfying the bound $\|V\|\leq d/2$. Moreover, the operator $L$
turns out to be similar to a self-adjoint operator whenever the
strict inequality \eqref{Vd2I} holds. These results combined in
Theorem~\ref{Tpi}\;(ii) below (see also Remark \ref{d2sharp})
represent an extension of similar results previously known due to
\cite{AdL1995} and \cite{MenShk} for the spectral dispositions with
subordinated $\sigma_0$ and $\sigma_1$.

By using the results of \cite{LT2004-PT,TretterWag2003}, we give a
positive answer to the question (i) also in the generic case (G)
provided that the unperturbed operator $A$ is bounded (see Theorem
\ref{Tspreal}). For $A$ unbounded, we prove that in case (G) the
spectrum of $L=A+V$ for sure is purely real if $V$ satisfies a
stronger bound $\|V\|\leq d/\pi$. The strict bound $\|V\|<d/\pi$
guarantees, in addition, that $L$ is similar to a self-adjoint
operator $\bigl($see Theorem~\ref{Tpi}\;(i)$\bigr)$. The question
whether this is true for $d/\pi\leq\|V\|<d/2$ remains an open
problem.

We answer the question (iii) by using the concept of the operator
angle between two subspaces (for discussion of this notion and
references see, e.g., \cite{KMM2}). Recall that if $\fM$ and $\fN$
are subspaces of a Hilbert space, the operator angle
$\Theta(\fM,\fN)$ between $\fM$ and $\fN$ measured relative to the
subspace $\fM$ is introduced by the following formula \cite{KMM3}:
\begin{equation}
\label{Theta}
\Theta(\fM,\fN)=\arcsin\sqrt{I_\fM-P_\fM P_\fN\bigl|_\fM},
\end{equation}
where $I_\fM$ denotes the identity operator on $\fM$ and
$P_\fM$ and $P_\fN$ stand for the orthogonal projections onto
$\fM$ and $\fN$, respectively.

Set
\begin{equation}
\label{dcases}
\delta=\left\{\begin{array}{cl}
2d/\pi, & \text{case (G)},\\
 d, & \text{case (S)},
\end{array}\right.
\end{equation}
and assume that $\|V\|<\delta/2$. Since in both the cases (G) and
(S) under this assumption we have got the positive answer to the
question (ii), one can easily identify the spectral subspaces
$\fH'_0$ and $\fH'_1$ of $L$ associated with the corresponding
perturbed spectral sets $\sigma'_0$ and $\sigma'_1$ (cf. Lemma
\ref{Tcomp}). Let $\Theta_j=\Theta(\fH_j,\fH'_j)$, $j=0,1,$ be the
operator angle between the unperturbed spectral subspace $\fH_j$ and
the perturbed one, $\fH'_j$. Our main result (presented in Theorem
\ref{Tpi}) regarding the operator angles $\Theta_0$ and $\Theta_1$
is that under condition $\|V\|<\delta/2$ the following bound holds:
\begin{equation}
\label{ThbpIn} \tan\Theta_j\leq\tanh\left(\frac{1}{2}\mathop{\rm
arctanh}\frac{2\|V\|}{\delta}\right),\quad j=0,1,
\end{equation}
which means, in particular, that $\Theta_j<\frac{\pi}{4}$, $j=0,1$.
Theorem \ref{Tpi} also gives the bounds on location of the perturbed
spectral sets $\sigma'_0$ and $\sigma'_1$ $\bigl($see formulas
\eqref{sigp}$\bigr)$.

In the case (S) the bounds on $\sigma'_0$ and $\sigma'_1$ as well as
the bounds \eqref{ThbpIn} are optimal (see Remark~\ref{Rsharp}).
Inequalities \eqref{ThbpIn} resemble the sharp norm estimate for the
operator angle between perturbed and unperturbed spectral subspaces
from the celebrated Davis-Kahan $\tan2\Theta$ Theorem (see
\cite{DK70}, p. 11; cf. \cite[Theorem 2.4]{KMM5} and \cite[Theorem
1]{MotSel}). Recall that the latter theorem serves for the case
where the unperturbed spectral subsets $\sigma_0$ and $\sigma_1$ are
subordinated and the off-diagonal perturbation $V$ is self-adjoint.
The difference is that the usual tangent of the Davis-Kahan
$\tan2\Theta$ Theorem is replaced on the right-hand side of
\eqref{ThbpIn} by the hyperbolic one. Another distinction is that
the bound \eqref{ThbpIn} holds not only for the subordinated
spectral sets $\sigma_0$ and $\sigma_1$ but also for the disposition
where one of these sets lies in a finite gap of the other set and thus
$\sigma_0$ and $\sigma_1$ are not subordinated.

The results obtained are of particular interest for the theory of
operators on Krein spaces \cite{AI}. The reason for this is that
introducing an indefinite inner product $[x,y]=(Jx,y)$, $x,y\in\fH$,
instead of the initial inner product $(\cdot,\cdot)$, turns $\fH$
into a Krein space. The operators $V$ and $L=A+V$ being
$J$-self-adjoint on $\fH$ appear to be self-adjoint operators on the
newly introduced Krein space $\fK$. Under the condition
$\|V\|<\delta/2$ in both cases (G) and (S) we establish that the
perturbed spectral subspaces $\fH'_0$ and $\fH'_1$ are mutually
orthogonal with respect to the inner product $[\cdot,\cdot]$.
Moreover, these subspaces are maximal uniformly positive and maximal
uniformly negative, respectively (see Remark \ref{RHHKr}). The
restrictions of $L$ onto $\fH'_0$ and $\fH'_1$ are $\fK$-unitary
equivalent to self-adjoint operators on $\fH_0$ and $\fH_1$,
respectively. This extends similar results previously known from
\cite{AdL1995} and \cite{MenShk} for the case where the spectral
sets $\sigma_0$ and $\sigma_1$ are subordinated.

Another motivation for the present paper is in the spectral analysis
of non-self-adjoint Schr\"o\-din\-ger operators that involve the
so-called $\cP\cT$-symmetric potentials. Starting from the
pioneering works \cite{Bender1,Bender2}, these potentials attracted
considerable attention because of their property to produce, in some
cases, purely real spectra (see, e.g.,
\cite{AlbFK,AlbKuzh2004,BenderRep,Krej,Mostaf,Znoj}). The local
$\cP\cT$-symmetric potentials appear to be $J$-self-adjoint with
respect to the space parity operator $\cP$ (see, e.g.,
\cite{LT2004-PT,Mostaf}), allowing for an embedding the problem into
the context of the spectral theory for $J$-self-adjoint
perturbations (this also means that the $\cP\cT$-symmetric
perturbations may be studied within the framework of the Krein space
theory \cite{AlbKuzh2004,LT2004-PT,Tanaka2}).

The main tool we use in our analysis is a reduction of the problems
(i)--(iii) to the study of the operator Riccati equation
$$
KA_0+A_1K+KBK=-B^*
$$
associated with the representation of the perturbed operator $L=A+V$
in the $2\times 2$ block matrix form
\begin{equation*}
%
L=\left(\begin{array}{rl}
  A_0      &  B \\
  -B^*     &   A_1
\end{array}\right),
\end{equation*}
where $A_0=A\bigr|_{\fH_0}$, $A_1=A\bigr|_{\fH_1}$, and
$B=V\bigr|_{\fH_1}$. Assuming \eqref{dcases}, we prove that the
Riccati equation has a bounded solution $K$ for any $B$ such
that $\|B\|<\delta/2$. The key statement is that the perturbed
spectral subspaces $\fH'_0$ and $\fH'_1$ are the graphs of the
operators $K$ and $K^*$, respectively, which then allows us
to derive the bounds \eqref{ThbpIn}.

The plan of the paper is as follows. In Section \ref{SecOR} we give
necessary definitions and present some basic results on the operator
Riccati equations associated with a class of unbounded
non-self-adjoint  $2\times 2$ block operator matrices. Section
\ref{SecSO} is devoted to the related Sylvester equations. In
Section \ref{SecROs} we prove a number of existence and uniqueness
results for the operator Riccati equations. In Section \ref{SecJsym}
we consider $J$-self-adjoint perturbations and find conditions on
their norm guaranteeing the reality of the resulting spectrum. In
this section we also prove the bound \eqref{ThbpIn} on the variation
of the spectral subspaces and discuss the embedding of the problem
into the context of the Krein space theory. Finally, in Section
\ref{SecExHO} we apply some of the results obtained to a
quantum-mechanical Hamiltonian describing the harmonic oscillator
under a $\mathcal{PT}$-symmetric perturbation.

We conclude the introduction with the description of some more
notations that are used thro\-ug\-h\-o\-ut the paper. By a subspace
we always understand a closed linear subset of a Hilbert space. The
identity operator on a subspace (or on the whole Hilbert space)
$\fM$ is denoted by $I_\fM$. If no confusion arises, the index $\fM$
may be omitted in this notation.  The Banach space of bounded linear
operators from a Hilbert space $\fM$ to a Hilbert space $\fN$ is
denoted by $\cB(\fM,\fN)$. For $\cB(\fM,\fM)$ we use a shortened
notation $\cB(\fM)$. By $\fM\oplus\fN$ we will understand the
orthogonal sum of two Hilbert spaces (or orthogonal subspaces) $\fM$
and $\fN$. By $\cO_r(\fM,\fN)$, $0\leq r<\infty$, we denote the
closed ball in $\cB(\fM,\fN)$, having radius $r$ and being centered
at zero, that is,
$$
\cO_r(\fM,\fN)=\{K\in\cB(\fM,\fN)\,\,\big|\,\,\, \|K\|\leq r\}.
$$
If it so happens that $r=+\infty$, by $\cO_\infty(\fM,\fN)$ we will
understand the whole space $\cB(\fM,\fN)$. The notation $\mathop{\rm
conv}(\sigma)$ is used for the convex hull of a Borel set
$\sigma\subset\bbR$. By $O_r(\Omega)$, $r\geq 0$, we denote the closed
$r$-neighborhood of a Borel set $\Omega$ in the complex plane $\bbC$, i.e.
$O_r(\Omega)=\{z\in\bbC\big|\,\dist(z,\Omega)\leq r\}$.

\section{Operator Riccati equation}
\label{SecOR}

We start by recalling the concepts of weak, strong, and operator
solutions to the operator Ric\-ca\-ti equation (see \cite{AMM,AM01}).

\begin{definition}
\label{DefSolRic} Assume that  $A_0$ and $A_1$ are possibly unbounded
densely defined closed operators on the Hilbert spaces $\fH_0$ and
$\fH_1$, respectively. Let $B$ and $C$ be bounded operators from $\fH_1$
to $\fH_0$ and from $\fH_0$ to $\fH_1$, respectively.

A bounded operator $K\in\cB(\fH_0,\fH_1)$ is said to be a weak
solution of the Riccati equation
\begin{equation}
\label{RicABCD}
KA_0-A_1K+KBK=C
\end{equation}
if
$$
%
\begin{array}{c}
( KA_0x,y)-( Kx,A_1^*y)+( KBKx,y)=( Cx,y)\\[0.5em]
 \text{ for all } x\in \dom (A_0)
\text{ and } y\in \dom(A_1^*).
\end{array}
$$

A bounded operator $K\in\cB(\fH_0,\fH_1)$ is called a strong
solution of the Riccati equation \eqref{RicABCD} if
\begin{equation}
\label{ranric}
\ran\bigl({K}|_{\dom(A_0)}\bigr)\subset\dom(A_1)
\end{equation}
and
\begin{equation}
\label{rics}
KA_0x-A_1Kx+KBKx=Cx  \text{ \, for all \, } x\in \dom(A_0).
\end{equation}

Finally,  $K\in \cB(\fH_0,\fH_1)$ is said to be
an operator solution of the Riccati equation \eqref{RicABCD} if
\begin{equation}
\label{ranric1}
\ran(K)\subset \dom(A_1),
\end{equation}
the operator $KA_0$ is bounded on $\dom(KA_0)=\dom(A_0)$,
and the equality
\begin{equation}
\label{Ricext}
\overline{KA_0}-A_1K+KBK=C
\end{equation}
holds as an operator equality, where $\overline{KA_0}$ denotes the closure of
$KA_0$.
\end{definition}
\begin{remark}
\label{RemRAdj}
We will call the equation
\begin{equation}
\label{RicAdj}
XA^*_1-A^*_0X-XB^*X=-C^*
\end{equation}
the adjoint of the operator Riccati equation \eqref{RicABCD}. It
immediately follows from the definition that an operator
$K\in\cB(\fH_0,\fH_1)$ is a weak solution to the Riccati equation
\eqref{RicABCD} if and only if the adjoint of $K$, $X=K^*$, is a
weak solution to the adjoint equation \eqref{RicAdj}.
\end{remark}

Clearly, any operator solution $K\in\cB(\fH_0,\fH_1)$ to the Riccati
equation \eqref{RicABCD} is automatically a strong solution.
Similarly, any strong solution is also a weak solution. But, in
fact, by a result of \cite{AM01} one does not need to distinguish
between weak and strong solutions to the Riccati equation
\eqref{RicABCD}. This is seen from the following statement.

\begin{lemma}[\cite{AM01}, Lemma 5.2]
\label{weak-strongR} Let $A_0$ and $A_1$ be densely defined possibly
unbounded closed operators on the Hilbert spaces $\fH_0$ and $\fH_1$,
respectively, and $B\in\cB(\fH_1,\fH_0)$, $C\in\cB(\fH_0,\fH_1)$. If
$K\in\cB(\fH_0,\fH_1)$ is a weak solution of the Riccati equation
\eqref{RicABCD} then $K$ is also a strong solution of
\eqref{RicABCD}.
\end{lemma}

If the operators $A_0$, $A_1$, $B$, and $C$ are as in Definition
\ref{DefSolRic} then a $2\times2$ operator block matrix
\begin{equation}
\label{L}
L=\left(\begin{array}{ll}
  A_0      &  B \\
  C   &   A_1
\end{array}\right), \qquad \dom(L)=\dom(A_0)\oplus\dom(A_1),
\end{equation}
is a densely defined and possibly unbounded closed operator on the
Hilbert space
\begin{equation}
\label{fHsum}
\fH=\fH_0\oplus\fH_1.
\end{equation}
The operator $L$ will often be viewed as the result of the
perturbation of the block diagonal matrix
\begin{equation}
\label{Adiag}
A=\mathop{\rm diag}(A_0,A_1),\qquad \dom(A)=\dom(A_0)\oplus\dom(A_1),
\end{equation}
by the off-diagonal bounded perturbation
\begin{equation}
\label{Voff}
V=\left(\begin{array}{ll}
  0   &  B \\
  C   &   0
\end{array}\right).
\end{equation}

The operator Riccati equation
\eqref{RicABCD} and the block operator matrix $L$ are usually said
to be associated to each other. Surely, one can also associate with
the matrix $L$ another operator Riccati equation,
\begin{equation}
\label{RicABCD1}
K'A_1-A_0K'+K'CK'=B,
\end{equation}
assuming that a solution $K'$ (if it exists) should be a bounded operator
from $\fH_1$ to $\fH_0$.

It is well known that the solutions to the Riccati equations \eqref{RicABCD}
and \eqref{RicABCD1} determine invariant subspaces for the operator matrix
$L$ (see, e.g., \cite{AMM} for the case where the matrix $L$ is self-adjoint
or \cite{LMMT} for the case of a non-self-adjoint $L$). These subspaces
have the form of the graphs
\begin{equation}
\label{cGK}
\cG(K)=\{x\in\fH_0\oplus\fH_1\,|\,x=x_0\oplus Kx_0 \text{\, for some
\,} x_0\in\fH_0\}
\end{equation}
and
\begin{equation}
\label{cGK1}
\cG(K')=\{x\in\fH_0\oplus\fH_1\,|\,x=K'x_1\oplus x_1 \text{\, for some
\,} x_1\in\fH_1\}
\end{equation}
of the corresponding (bounded) solutions $K$ and $K'$. Notice that
the subspaces of the form \eqref{cGK} and \eqref{cGK1} are usually
called the graph subspaces associated with the operators $K$ and
$K'$, respectively, while $K$ and $K'$ themselves are called the
angular operators. Usage of the latter term is explained, in
particular, by the fact that if a subspace $\fG\subset\fH$ is a
graph $\fG=\cG(K)$ of a bounded linear operator $K$ from a subspace
$\fM$ to its orthogonal complement $\fM^\perp$,
$\fM^\perp=\fH\ominus\fM$, then the following equality holds (see
\cite{KMM2}; cf. \cite{DK70} and \cite{Halmos:69}):
\begin{equation}
\label{KThet}
|K|=\tan\Theta(\fM,\fG),
\end{equation}
where $|K|$ is the absolute value of $K$, $|K|=\sqrt{K^*K}$,  and
$\Theta(\fM,\fG)$ the operator angle \eqref{Theta} between
the subspaces $\fM$ and $\fG$.

The precise statement relating solutions of the Riccati
equations \eqref{RicABCD} and \eqref{RicABCD1} to invariant
subspaces of the operator matrix \eqref{L} is as follows.

\begin{lemma}
\label{Lgraph} Let the entries $A_0$, $A_1$, $B$, and $C$ be as in
Definition \ref{DefSolRic} and let a $2\times2$ block operator
matrix $L$ be given by \eqref{L}. Then the graph $\cG(K)$ of a
bounded operator $K$ from $\fH_0$ to $\fH_1$ satisfying
\eqref{ranric} is an invariant subspace for the operator matrix $L$
if and only if $K$ is a strong solution to the operator Riccati
equation \eqref{RicABCD}. Similarly, the graph $\cG(K')$ of an
operator $K'\in\cB(\fH_1,\fH_0)$ such that
$\ran\bigl({K'}|_{\dom(A_1)}\bigr)\subset\dom(A_0)$ is an invariant
subspace for $L$ if and only if this operator is a strong solution
to the Riccati equation \eqref{RicABCD1}.
\end{lemma}

The proof of this lemma is straightforward and follows the same line as
the proof of the corresponding part in \cite[Lemma 5.3]{AMM}. Thus, we omit it.

The next assertion contains two useful identities involving the
strong solutions to the Riccati equations \eqref{RicABCD} and \eqref{RicABCD1}.

\begin{lemma}
\label{L2ident} Let the entries $A_0$, $A_1$, $B$, and $C$ be as in
Definition \ref{DefSolRic}. Assume that operators
$K\in\cB(\fH_0,\fH_1)$ and $K'\in\cB(\fH_1,\fH_0)$ are strong
solutions to equations \eqref{RicABCD} and \eqref{RicABCD1},
respectively. Then
\begin{equation}
\label{DomKK}
\ran\bigl({K'K}|_{\dom(A_0)}\bigr)\subset\dom(A_0), \quad
\ran\bigl({KK'}|_{\dom(A_1)}\bigr)\subset\dom(A_1),
\end{equation}
and
\begin{align}
\label{Id0}
(I-K'K)(A_0+BK)x &=(A_0-K'C)(I-K'K)x \quad\text{for all \,}x\in\dom(A_0),\\
\label{Id1}
(I-KK')(A_1+CK')y &=(A_1-KB)(I-KK')y \quad\text{for all \,}y\in\dom(A_1).
\end{align}
\end{lemma}
\begin{proof}
The inclusions \eqref{DomKK} follow immediately from the definition
of a strong solution to the operator Riccati equation (see condition
\eqref{ranric}).

Let $x\in\Dom(A_0)$. Taking into account the first of the inclusions
\eqref{DomKK} as well as the inclusions
$\ran\bigl({K}|_{\dom(A_0)}\bigr)\subset\dom(A_1)$ and
$\ran\bigl({K'}|_{\dom(A_1)}\bigr)\subset\dom(A_0)$ one can write
\begin{align}
\nonumber
(A_0-K'C)(I-K'K)x &=(A_0-K'C)x-(A_0K'-K'CK')Kx\\
\label{Id0a}
&=(A_0-K'C)x-(K'A_1-B)Kx,
\end{align}
by making use of the Riccati equation \eqref{RicABCD1} itself at the second step.
Similarly,
\begin{align}
\nonumber
(I-K'K)(A_0+BK)x &=(A_0+BK)x-K'(KA_0+KBK)x\\
\label{Id0b}
&=(A_0+BK)x-K'(C+A_1K)x,
\end{align}
due to the Riccati equation \eqref{RicABCD}. Comparing \eqref{Id0a} and \eqref{Id0b}
we arrive at the identity \eqref{Id0}.

Identity \eqref{Id1} is proven analogously.
\end{proof}

We will also need the following auxiliary lemma.
\begin{lemma}
\label{Lsum} Suppose that the operators $K\in\cB(\fH_0,\fH_1)$ and
$K'\in\cB(\fH_1,\fH_0)$ are such that the $2\times2$ operator block
matrix
\begin{equation}
\label{W}
W=\left(\begin{array}{ll}
  I      &  K' \\
  K   &   I
\end{array}\right)
\end{equation}
considered on $\fH=\fH_0\oplus\fH_1$ is boundedly invertible, i.e.
the inverse operator $W^{-1}$ exists and is bounded.
Then the graphs $\cG(K)$ and $\cG(K')$ of the
operators $K$ and $K'$ are linearly independent subspaces of $\fH$
and
\begin{equation}
\label{HKK}
\fH=\cG(K)\dotplus\cG(K'),
\end{equation}
where the sign ``$\dotplus$'' denotes the direct sum of two subspaces.
\end{lemma}
\begin{proof}
The existence and boundedness of $W^{-1}$ imply that equation $Wx=y$
is uniquely solvable for any $y\in\fH$. This means that there are
unique $x_0\in\fH_0$ and unique $x_1\in\fH_1$ such that $y=x_0\oplus
Kx_0+K'x_1\oplus x_1$ and hence $\fH\subset\cG(K)+\cG(K')$. Since
both $\cG(K)$ and $\cG(K')$ are subspaces of $\fH$, the inclusion
turns into equality, $\fH=\cG(K)+\cG(K')$. The linear independence
of  $\cG(K)$ and $\cG(K')$ follows from the fact that equation
$Wx=0$ has only the trivial solution $x=0$.
\end{proof}

\begin{remark}
\label{RemKK}
It is well known that the following three statements are equivalent.
\begin{enumerate}
\item[(i)] The operator matrix \eqref{W} is boundedly invertible.
\item[(ii)] The inverse $(I-KK')^{-1}$ exists and  is bounded.
\item[(iii)] The inverse $(I-K'K)^{-1}$ exists and is bounded
\end{enumerate}
For a proof of this assertion see, e.g. \cite[Theorem 1.1. and Lemma
2.1]{HKM} where even a Banach-space case of $2\times2$ block
operator matrices of the form \eqref{W} with unbounded entries $K$
and $K'$ has been studied.
\end{remark}

\begin{remark}
The inverse of the operator $W$ is explicitly written as
\begin{align}
\nonumber
W^{-1}&=\begin{pmatrix}
(I-K'K)^{-1} & -K'(I-KK')^{-1}\\
-K(I-K'K)^{-1} & (I-KK')^{-1}
\end{pmatrix}\\
\label{W-1a}
&=
\begin{pmatrix}
(I-K'K)^{-1} & -(I-K'K)^{-1}K'\\
-(I-KK')^{-1}K & (I-KK')^{-1}.
\end{pmatrix}
\end{align}
The (oblique) projections $Q_{\cG(K)}$ and $Q_{\cG(K')}$ onto the graph subspaces
$\cG(K)$ and $\cG(K')$ along the corresponding complementary graph subspaces $\cG(K')$
and $\cG(K)$ are given by
\begin{equation}
\label{PGK}
Q_{\cG(K)}=\begin{pmatrix} I \\ K \end{pmatrix}
(I-K'K)^{-1}
\begin{pmatrix}
I & -K'
\end{pmatrix}
\text{\, and \,}
Q_{\cG(K')}=\begin{pmatrix} K' \\ I \end{pmatrix}
(I-KK')^{-1}
\begin{pmatrix}
-K & I
\end{pmatrix},
\end{equation}
respectively.
\end{remark}

\begin{corollary}
\label{Ldiag} Assume the hypothesis of Lemma \ref{Lgraph}. Suppose that
$K\in\cB(\fH_0,\fH_1)$ and $K'\in\cB(\fH_1,\fH_0)$ are strong
solutions to the Riccati equations \eqref{RicABCD} and
\eqref{RicABCD1}, respectively. Assume, in addition, that the
$2\times2$ operator block matrix $W$ formed with these solutions
according to \eqref{W} is a boundedly invertible operator on
$\fH=\fH_0\oplus\fH_1$. Then:
\begin{enumerate}
\item[(i)] The operator $L$ is
similar to a block diagonal operator matrix $Z=\diag(Z_0,Z_1)$,
\begin{equation}
\label{LWZ}
L=WZW^{-1},
\end{equation}
where $Z_0$ and $Z_1$ are operators on $\fH_0$ and $\fH_1$, respectively,
given by
\begin{align}
\label{Z0}
Z_0&=A_0+BK, \quad \dom(Z_0)=\dom(A_0),\\
\label{Z1}
Z_1&=A_1+CK', \quad \dom(Z_1)=\dom(A_1).
\end{align}
\item[(ii)] The Hilbert space $\fH$ splits into the
direct sum  $\fH=\fH'_0\dotplus\fH'_1$ of the graph subspaces
$\fH'_0=\cG(K)$ and $\fH'_1=\cG(K')$ that are invariant under $L$.
The restrictions $L|_{\fH'_0}$ and $L|_{\fH'_1}$ of $L$ onto
$\fH'_0$ and $\fH'_1$ are similar to the operators $Z_0$ and $Z_1$,
\begin{equation}
\label{W01Z}
W_0^{-1}L|_{\fH'_0}W_0=Z_0\quad\text{and}\quad W_1^{-1}L|_{\fH'_1}W_1=Z_1,
\end{equation}
where the entries $W_0:\,\fH_0\to\fH'_0$ and $W_1:\,\fH_1\to\fH'_1$
correspond to the respective columns of the block operator matrix
$W$,
\begin{equation}
\label{W01}
W_0 x_0=\begin{pmatrix} I \\ K\end{pmatrix}x_0,
\,\, x_0\in\fH_0,\quad\text{and}\quad
W_1 x_1=\begin{pmatrix} K' \\ I\end{pmatrix}x_1,
\,\, x_1\in\fH_1.
\end{equation}
\end{enumerate}
\end{corollary}
\begin{proof}
First, one verifies by inspection that $LW=WZ$ taking to account
that $K$ and $K'$ are the strong solutions to the Riccati equations
\eqref{RicABCD} and \eqref{RicABCD1}, respectively. The remaining statements
immediately follow from Lemma \ref{Lgraph} combined with Lemma \ref{Lsum}.
\end{proof}
\begin{remark}
\label{RZZ}
The similarity \eqref{LWZ} of the operators $L$ and $Z$ implies that
the spectrum of $L$ coincides with the union of the spectra of $Z_0$
and $Z_1$, that is, $\spec(L)=\spec(Z_0)\cup\spec(Z_1)$.
\end{remark}

\section{Operator Sylvester equation}
\label{SecSO}

Along with the Riccati equation \eqref{RicABCD}
we need to consider the operator Sylvester equation
\begin{equation}
\label{Syl}
XA_0-A_1X=Y
\end{equation}
assuming that the entries $A_0$ and $A_1$ are as in Definition
\ref{DefSolRic} and $Y\in\cB(\fH_0,\fH_1)$. The Sylvester equation
is a particular (linear) case of the Riccati equation and its weak,
strong, and operator solutions $X\in\cB(\fH_0,\fH_1)$ are understood
in the same way as in the above definition. Furthermore, by Lemma
\ref{weak-strongR} (cf. \cite[Lemma 1.3]{ARS1994}) one does not need
to distinguish between the weak and strong solutions to \eqref{Syl}.

Because of its importance for various areas of mathematics there is
an enormous literature on the Sylvester equation (for a review and
many references see paper \cite{BR}). With equation \eqref{Syl} one
often associates the Sylvester operator $\sS$ defined on the Banach space
$\cB(\fH_0,\fH_1)$ by the left-hand side of  \eqref{Syl}:
\begin{equation}
\label{SylOp}
\sS(X)=XA_0-A_1X
\end{equation}
with domain
\begin{equation}
\label{SylDom}
\Dom(\sS)=\left\{X\in\cB(\fH_0,\fH_1)\,\,\big|\,\,
\Ran(X\big|_{\Dom(A_0)})\subset\Dom(A_1)\right\}.
\end{equation}
Clearly, the Sylvester equation \eqref{Syl} has a unique solution
$X\in\Dom(\sS)$ if and only if $0\not\in\spec(\sS)$. It is known
that in general the spectrum of $\sS$ is larger than the (numerical)
difference between the spectra of $A_0$ and $A_1$. More precisely,
provided that $\spec(A_0)\neq\bbC$ or $\spec(A_1)\neq\bbC$ always
the following inclusion holds \cite{ARS1994}:
\begin{equation}
\label{sTAA}
\overline{\spec(A_0)-\spec(A_1)}\subset\spec(\sS),
\end{equation}
where we use the notation $\Sigma-\Delta=\left\{z-\zeta\,\,|\,\,
z\in\Sigma,\zeta\in\Delta \right\}$ for the numerical difference
between two Borel subsets $\Sigma$ and $\Delta$ of the
complex plane $\bbC$. The opposite inclusion in \eqref{sTAA}
may fail to hold if both operators $A_0$ and $A_1$ are unbounded.
The corresponding example was first given by V.\,Q.\,Ph\'{o}ng \cite{P91}
for the Sylvester equation \eqref{Syl} where one of the entries
$A_0$ and $A_1$ is an operator on a Banach (but not Hilbert) space.
An example where both $A_0$ and $A_1$ are operators on Hilbert
spaces and $\spec(\sS)\not\subset\overline{\spec(A_0)-\spec(A_1)}$
may be found in \cite[Example 6.2]{ARS1994}. Equality
\begin{equation}
\label{SAAeq}
\spec(\sS)=\spec(A_0)-\spec(A_1)
\end{equation}
holds if both $A_0$ and $A_1$ are bounded operators. This result is
due to G.~Lumer and M.~Rosenblum \cite{LumerR}. Equality \eqref{SAAeq}
also holds if only one of the entries $A_0$ and $A_1$ is a bounded
operator \cite{ARS1994}. In this case \eqref{SAAeq}
implies that if the spectra $A_0$ and $A_1$
are disjoint then $0\not\in\spec(\sS)$ and hence the operator $\sS$
is boundedly invertible. Moreover, a unique solution of the Sylvester
equation \eqref{Syl} admits an ``explicit'' representation
in the form a contour integral.

\begin{lemma}
\label{Krein}
Let $A_0$ be a possibly unbounded densely defined closed operator
on the Hilbert space $\fH_0$ and $A_1$ a  bounded  operator on the
Hilbert space $\fH_1$ such that
$$
\spec(A_0)\cap\spec(A_1)=\emptyset
$$
and $Y\in \cB(\fH_0,\fH_1)$.  Then  the Sylvester equation
\eqref{Syl} has a unique operator solution
\begin{equation}
\label{ESylKrein}
X=\frac{1}{2\pi\ri}
\int_{\gamma} dz\,(A_1-z)^{-1}Y (A_0-z)^{-1},
\end{equation}
where $\gamma$ is a union of closed contours in $\bbC$ with total
winding numbers $0$ around $\spec(A_0)$ and $1$ around $\spec(A_1)$ and
the integral converges in the norm operator topology.
\end{lemma}
\begin{corollary}
\label{CorKrein}
Under the hypothesis of Lemma \ref{Krein} the norm of the
inverse of the Sylvester operator $\sS$ may be estimated as
$$
\|\sS^{-1}\|\le (2\pi)^{-1}|\gamma|
\sup_{z\in \gamma} \|(A_0-z)^{-1}\|\, \|(A_1-z)^{-1}\|\,,
$$
where $|\gamma|$ denotes the length of the contour $\gamma$ in
\eqref{ESylKrein}.
\end{corollary}

The result of Lemma \ref{Krein} may be attributed to M.\,G.\,Krein
who lectured on the operator Sylvester equation in late 1940s  (see
\cite{BR}).  Later, it was independently obtained by
M.\,Ro\-sen\-blum~\cite{R56}.

As for the Sylvester operator \eqref{SylOp} with both unbounded
entries $A_0$ and $A_1$, we have an important result which is due
to W. Arendt, F. R\"abiger, and A. Sourour (see \cite[Theorem~4.1
and Corollary 5.4]{ARS1994}).

\begin{theorem}[\cite{ARS1994}]
\label{Arendt} Let $A_0$ and $A_1$ be closed densely defined
operators on the Hilbert spaces $\fH_0$ and $\fH_1$, respectively.
Assume that one (or both) of the following holds (hold) true:
\begin{enumerate}
\item[(i)] $A_0$ and $(-A_1)$ are generators of eventually norm continuous
$C_0$-semigroups;

\item[(ii)] $A_0$ and $(-A_1)$ are generators of $C_0$-semigroups
one of which is holomorphic.
\end{enumerate}
Then the spectrum of the Sylvester operator \eqref{SylOp} is
given by \eqref{SAAeq}.

\end{theorem}

Recall that an operator $H$ on the Hilbert space $\fM$ is said to be
$m$-dissipative if it is closed and both the spectrum and numerical
range of $H$ are contained in the left half-plane
\mbox{$\{z\in\bbC\,|$} \mbox{$\Img z\leq 0\}$}. The Lumer-Phillips
theorem asserts (see, e.g., \cite[Section II.3.b]{EN2000}; cf.
\cite[Theorem~B.21]{Haase}) that a $C_0$-semigroup on $\fM$ is a
contraction semigroup if and only if its generator is an
$m$-dissipative operator. The next statement represents a
generalization of a well known result by E. Heinz
(\cite[Satz~5]{H51}) to the case of unbounded operators. Notice that
the exponential ${\rm e}^{H t}$, $t\geq 0$, is understood below as
the corresponding element of the strongly continuous contraction
semigroup generated by an (unbounded) $m$-dissipative operator $H$.

\begin{theorem}
\label{Kexp} Let $A_0+\frac{\delta}{2}I$ and
$-A_1+\frac{\delta}{2}I$, $\delta>0$, be $m$-dissipative operators
on the Hilbert spaces $\fH_0$ and $\fH_1$, respectively, and $Y\in
\cB(\fH_0,\fH_1)$.  Then  the Sylvester equation \eqref{Syl} has a
uni\-que weak (and hence unique strong) solution given by
\begin{equation}
\label{XYSol1}
 X=-\int_0^{+\infty}
dt\,{\rm e}^{-A_1 t}Y{\rm e}^{A_0 t},
%
%
\end{equation}
where the integral is understood in the weak operator topology.
Moreover, the norm of the so\-lu\-tion \eqref{XYSol1} satisfies the
estimate
\begin{equation}
\label{XdY}
\| X\|\leq \frac{1}{\delta}\,\|Y\|.
\end{equation}
\end{theorem}
\begin{proof}
Under the hypothesis the operators $A_0$ and $(-A_1)$ are
themselves $m$-dissipative. Let $U_0(t)$ and $U_1(t)$, $t\geq 0$, be
contraction $C_0$-semigroups generated respectively by $A_0$ and
$(-A_1)$, that is, $U_0(t)={\rm e}^{A_0 t}$ and $U_1(t)={\rm e}^{-A_1
t}$. Clearly,
\begin{equation}
\label{UU}
\|U_0(t)\|\leq{\rm e}^{-\frac{\delta}{2}t}
\text{\, and \,}\|U_1(t)\|\leq{\rm e}^{-\frac{\delta}{2}t},\quad t\geq 0.
\end{equation}
The same bound also holds for the adjoint semigroup $U_1(t)^*$,
$t\geq 0$, whose generator is the \mbox{$m$-dissipative} operator
$(-A_1^*)$. Pick up arbitrary $x\in\fH_0$ and $y\in\fH_1$ and
introduce the orbit maps $t\mapsto\xi_x(t)=U_0(t)x$ and $\zeta_y:
t\mapsto\zeta_y(t)=U_1(t)^*y$. By the definition of a strongly
continuous semigroup, these maps are continuous functions of
$t\in[0,\infty)$. Taking into account the bounds \eqref{UU} one then
concludes that the improper integral
\begin{equation}
\int_0^\infty dt\, (U_1(t)YU_0(t)x,y)=\int_0^\infty dt\,
(Y\xi_x(t),\zeta_y(t))
\end{equation}
converges and its absolute value is bounded by
$\|Y\|\|x\|\|y\|/\delta$. Thus, the weak integral on the right-hand side of
\eqref{XYSol1} exists and the bound \eqref{XdY} holds.

Now assume that $x\in\Dom(A_0)$ and $y\in\Dom(A_1^*)$. In this case
the orbit maps $\xi_x(t)$ and $\zeta_y(t)$ are continuously
differentiable in $t$ and $\frac{d}{dt}\xi_x(t)=U_0(t)A_0x$,\,
$\frac{d}{dt}\zeta_y(t)=-U_1(t)^*A_1^*y$. For $X$ given by
\eqref{XYSol1}, an elementary computation shows that
\begin{align*}
(XA_0x,y)-(Xx,A_1^*y)&=-\int_0^\infty dt\,
\bigl((Y\mbox{$\frac{d}{dt}$}\xi_x(t),\zeta_y(t))+
(Y\xi_x(t),\mbox{$\frac{d}{dt}\zeta_y(t)$})\bigr)\\
&=-\int_0^\infty dt\,\mbox{$\frac{d}{dt}$}
(Y\xi_x(t),\zeta_y(t))=(Y\xi(0),\zeta(0)),
\end{align*}
taking into account  \eqref{UU} in the last step. Since $\xi_x(0)=x$
and $\zeta_y(0)=y$, by Definition \ref{DefSolRic} this implies that
the integral \eqref{XYSol1} is a weak (and hence strong) solution to
the Sylvester equation~\eqref{Syl}.

To prove the uniqueness of the weak solution \eqref{XYSol1} it is
sufficient to show that the homogeneous Sylvester equation
$XA_0-A_1X=0$ has the only weak solution $X=0$.
For a weak solution $X$ to this equation we have
\begin{equation}
\label{SylHom}
(XA_0u,v)-(Xu,A_1^*v)=0\quad\text{for all \,}u\in\Dom(A_0)\text{\,
and \,}v\in\Dom(A_1^*).
\end{equation}
Take the vectors $u$ and $v$ of the form $u=\xi_x(t)={\rm e}^{A_0
t}x$, $v=\zeta_y(t)={\rm e}^{-A_1^* t}y$, $t\geq0$, where the orbit
maps $\xi_x(t)$ and $\zeta_y(t)$ correspond to some $x\in\dom(A_0)$
and $y\in\dom(A_1^*)$ and hence are both continuously
differentiable in $t\in[0,\infty)$. Notice that the assumption $x\in\dom(A_0)$,
$y\in\dom(A_1^*)$ also implies (see, e.g. \cite[Lemma 1.3]{EN2000})
that $u\in\Dom(A_0)$, $v\in\Dom(A_1)$, and
$A_0u=\frac{d}{dt}\xi_x(t)$, \mbox{$A_1^*v=-\frac{d}{dt}\zeta(t)$}.
With such a choice of $u$ and $v$ it follows from \eqref{SylHom}
that
\begin{equation*}
\mbox{$\frac{d}{dt}$}(X\xi_x(t),\zeta_y(t))=0\quad\text{whenever
\,}x\in\Dom(A_0)\text{ and }y\in\Dom(A_1^*).
\end{equation*}
Hence the function $(X\xi_x(t),\zeta_y(t))$, $t\geq 0$, is a
constant. Moreover, it equals zero since it vanishes
as~$t\to\infty$. This yields in particular that
\begin{equation*}
(Xx,y)=(X\xi_x(0),\zeta_y(0))=0 \quad\text{for all
\,}x\in\Dom(A_0)\text{ and }y\in\Dom(A_1^*).
\end{equation*}
The latter implies $X=0$, which completes the proof.
\end{proof}
\begin{remark}
A statement similar to Theorem \ref{Kexp} was previously announced
without a proof in~\cite{AMM} (see \cite[Lemma 2.6]{AMM}).
\end{remark}

The second important example where a bound of the  \eqref{XdY} type
exists is given in \cite[Theorem~3.2]{BDM1983}. This example is
as follows.
\begin{theorem}
\label{Th-BDM1983}
Assume that the operators $A_0$ and $A_1$ are densely defined and closed.
Assume, in addition, that there is $\lambda\in\varrho(A_1)$ such that
$\|A_0-\lambda\|\leq r$ and $\|(A_1-\lambda)^{-1}\|\leq (r+\delta)^{-1}$
for some $r\geq 0$ and $\delta>0$. Then for any $Y\in\cB(\fH_0,\fH_1)$
the unique strong solution $X$
to the Sylvester equation \eqref{Syl} admits the estimate
$\delta\|X\|\leq \|Y\|.$
\end{theorem}
If the operators $A_0$ and $A_1$ are normal then no reference point
$\lambda$ is needed and the result is stated in a more universal form
(see \cite[Theorem~3.2]{BDM1983}).
\begin{corollary}
\label{Xnormal}
Let both $A_0$ and $A_1$ be normal operators such that $\spec(A_0)$
is contained in a closed disk of radius $r$, $r\geq 0$, while
$\spec(A_1)$ is disjoint from the open disk (with the same
center) of radius $r+\delta$, $\delta>0$. Then for any
$Y\in\cB(\fH_0,\fH_1)$ the Sylvester equation \eqref{Syl} has a
unique strong solution $X$ and $\delta\|X\|\leq \|Y\|$.
\end{corollary}

The above two theorems and the corollary give examples where the bounded
inverse of the Sylvester operator $\sS$ exists and for
the norm of $\sS^{-1}$ the estimate
$\delta\|\sS^{-1}\|\leq 1$
holds with some $\delta>0$. Moreover, this estimate is universal in
the sense that it remains valid for any $A_0$ and $A_1$ satisfying
the corresponding hypotheses.

\section{Existence results for the Riccati equation}
\label{SecROs}

In this section we return to the operator Riccati equation \eqref{RicABCD} to prove
some sufficient conditions for its solvability. In their proof we will
rely just on the assumption that an estimate like \eqref{XdY} holds
for the solution of the corresponding Sylvester equation.

\begin{theorem}
\label{ExistGen} Let  $A_0$ and $A_1$ be possibly unbounded closed
densely defined operators on the Hilbert spaces $\fH_0$ and $\fH_1$,
respectively. Assume that the Sylvester operator $\sS$ defined on
$\cB(\fH_0,\fH_1)$ by \eqref{SylOp} and
\eqref{SylDom} is boundedly invertible {\rm(}that is,
$0\not\in\spec(\sS)${\rm)} and
\begin{equation}
\label{S1del}
\|\sS^{-1}\|\leq \frac{1}{\delta}
\end{equation}
for some $\delta>0$. Assume, in addition, the operators
$B\in\cB(\fH_1,\fH_0)$ and $C\in\cB(\fH_0,\fH_1)$ are
such that the following bound holds:
\begin{equation}
\label{BCdel}
\sqrt{\|B\|\|C\|}<\dfrac{\delta}{2}.
\end{equation}
Then the Riccati equation \eqref{RicABCD} has a unique strong solution
in the ball $\cO_{\delta/(2\|B\|)}(\fH_1,\fH_0)$.
The strong solution $K$ satisfies the estimate
\begin{equation}
\label{EstL}
\|K\|\leq \dfrac{\|C\|}{\frac{\delta}{2}+\sqrt{\frac{\delta^2}{4}-
\|B\|\,\|C\|}}.
\end{equation}
\end{theorem}
\begin{proof}
If $B=0$ then the assertion, including the estimate \eqref{EstL},
follows immediately from the hypothesis on the invertibility of
$\sS$ on $\cB(\fH_0,\fH_1)$ taking into account the bound
\eqref{S1del}.

Suppose that $B\neq 0$. In this case the proof is performed by
applying Banach's Fixed Point Theorem. First, we notice
that the bounded invertibility of $\sS$ on $\cB(\fH_0,\fH_1)$ allows us to
rewrite the Riccati equation \eqref{RicABCD} in the form
$$
K=F(K)
$$
where the mapping $F:\cB(\fH_0,\fH_1)\to\Dom(\sS)$ is given by
$$
F(K)=\sS^{-1}(C-KBK).
$$
By \eqref{S1del} we have
\begin{equation}
\label{EstF1}
\|F(K)\|\leq \frac{1}{\delta}(\|C\|+\|B\|\,\|K\|^2),
 \quad K\in\cB(\fH_0,\fH_1)
\end{equation}
and
\begin{equation}
\label{EstF2}
\|F(K_1)-F(K_2)\|\leq \frac{1}{\delta}
\|B\|\,(\|K_1\|+\|K_2\|)\, \|K_1-K_2\|, \quad K_1,K_2\in\cB(\fH_0,\fH_1).
\end{equation}
The bound \eqref{EstF1} implies that $F$ maps the ball
$\cO_r(\fH_0,\fH_1)$
into itself whenever
\begin{equation}
\label{eq1}
\|B\|\,r^2+\|C\|\leq r\delta.
\end{equation}
At the same time, from \eqref{EstF2} it follows that $F$ is a strict contraction
of the ball  $\cO_r(\fH_1,\fH_0)$ whenever
\begin{equation}
\label{eq2}
2\|B\|r<\delta.
\end{equation}
Solving inequalities \eqref{eq1} and \eqref{eq2} one concludes that
if the radius  $r$ of the ball $\cO_r(\fH_1,\fH_0)$ is within the bounds
\begin{equation}
\label{IneqR}
\dfrac{\|C\|}{\frac{\delta}{2}+\sqrt{\frac{\delta^2}{4}-
\|B\|\,\|C\|}} \leq r <
\frac{\delta}{2\|B\|},
\end{equation}
then $F$ is a strictly contractive mapping of the ball
$\cO_r(\fH_1,\fH_0)$ into itself. Applying Banach's Fixed Point
Theorem one then infers that equation \eqref{RicABCD} has a unique
solution within any ball $\cO_r(\fH_1,\fH_0)$ whenever the radius $r$ satisfies
\eqref{IneqR}. This means that the fixed point is the same for all
the radii satisfying \eqref{IneqR} and hence it belongs to the
smallest of the balls. This conclusion proves the bound
\eqref{EstL} and completes the whole proof.
\end{proof}
\begin{remark}
In \eqref{BCdel}--\eqref{EstL} one may set
$\delta=\|\sS^{-1}\|^{-1}$.
\end{remark}
\begin{remark}
\label{Bhyp} By using the hyperbolic tangent function and its inverse the
bound \eqref{EstL} (for \mbox{$B\neq 0$}) can be equivalently written in
the hypertrigonometric form
\begin{equation}
\label{EstLt} \|K\|\leq
\sqrt{\dfrac{\|C\|}{\|B\|}}\,\,\tanh\left(\dfrac{1}{2}
\mathop{\rm arctanh}\dfrac{2\sqrt{\|B\|\|C\|}}{\delta}\right).
\end{equation}
Notice that under condition \eqref{BCdel} we always have
$$
\tanh\left(\dfrac{1}{2}
\mathop{\rm arctanh}\dfrac{2\sqrt{\|B\|\|C\|}}{\delta}\right)<1.
$$
\end{remark}
\begin{remark}
Fixed-point based approaches to prove the solvability of the
operator Riccati equation with bounded
entries $A_0$ and $A_1$ have been used in many papers (see, e.g.,
\cite{AdLT}, \cite{Demmel1987}, \cite{Nair2001}, \cite{Stewart1973},
\cite{Stewart1971}). In the case where at least one of the entries
$A_0$ and $A_1$ is an unbounded self-adjoint or normal operator, a
fixed-point approach has been employed in \cite{AMM}, \cite{AM01},
\cite{MM99}, and \cite{MotRem}. Theorem \ref{ExistGen} represents an
extension of the fixed-point existence results obtained in
\cite[Theorem 3.5]{Stewart1971} and \cite[Theorem 3.1]{Nair2001} for
the Riccati equation \eqref{RicABCD} with both bounded $A_0$ and
$A_1$ to the case where the entries $A_0$ and $A_1$ are not
necessarily bounded.
\end{remark}

\begin{theorem}
\label{LdiagF} Assume the hypothesis of Theorem \ref{ExistGen}. Then
the block operator matrix $L$ defined by \eqref{L} is block
diagonalizable with respect to the direct sum decomposition
$\fH=\cG(K)\dotplus\cG(K')$ where $K$ is the unique strong solution
to the Riccati equation \eqref{RicABCD} within the operator ball
$\cO_{\delta/(2\|B\|)}(\fH_0,\fH_1)$ and $K'$ the unique strong
solution to the Riccati equation \eqref{RicABCD1} within the
operator ball $\cO_{\delta/(2\|C\|)}(\fH_1,\fH_0)$.
\end{theorem}
\begin{proof}
By Theorem \ref{ExistGen} for $K$ the estimate \eqref{EstL} holds.
By the same theorem for $K'$ we have
\begin{equation}
\label{EstLp}
\|K'\|\leq {\|B\|}\left({\frac{\delta}{2}+\sqrt{\frac{\delta^2}{4}-
\|B\|\,\|C\|}}\right)^{-1}.
\end{equation}
Then the hypothesis $\|B\|\|C\|<\delta/2$ also implies that
$\|K\|\|K'\|<1$. Hence by Remark \ref{RemKK} the operator $W$ in
\eqref{W} is boundedly invertible. Applying Corollary \ref{Ldiag}
completes the proof.
\end{proof}

\begin{remark}
\label{Rrb}
If, in addition, both operators $A_0$ and $A_1$ are normal then,
for $r$ defined by
\begin{equation}
\label{rb}
r=\dfrac{\|B\|\,\|C\|}{\frac{\delta}{2}+\sqrt{\frac{\delta^2}{4}-
\|B\|\,\|C\|}}=\sqrt{\|B\|\,\|C\|}\,\,\tanh\left(\dfrac{1}{2}
\mathop{\rm arctanh}\dfrac{2\sqrt{\|B\|\|C\|}}{\delta}\right),
\end{equation}
the spectrum of the block matrix $L$ lies in the closed
$r$-neighborhood of the spectrum of its main-diagonal part
$A=\diag(A_0, A_1)$. That is,
$\dist\bigl(z,\spec(A_0)\cup\spec(A_1)\bigr)\leq r$ whenever
$z\in\spec(L)$. This immediately follows from the representation
\eqref{LWZ}--\eqref{Z1} and the bounds \eqref{EstL} and \eqref{EstLp}
(see also Remark \ref{RZZ}). Notice that if $B\neq 0$ and $C\neq 0$
then $r<\sqrt{\|B\|\|C\|}$ and hence $r<\|V\|$ taking into account
that $\|V\|=\max(\|B\|,\|C\|)$.
\end{remark}

{}From now on we assume that the entries $A_0$ and $A_1$ are self-adjoint
operators with disjoint spectra and thus adopt the following

\begin{hypothesis}
\label{HLss} Let  $A_0$ and $A_1$ be (possibly unbounded)
self-adjoint operators on the Hilbert spaces $\fH_0$ and $\fH_1$
with domains $\dom(A_1)$ and $\dom(A_1)$, respectively. Assume that
the spectra of the operators $A_0$ and $A_1$ are disjoint and let
\begin{equation}
\label{AAdist}
d=\dist\bigl(\spec(A_0),\spec(A_1)\bigr)\, \bigl(>0\bigr).
\end{equation}
\end{hypothesis}

Hypothesis \ref{HLss} imposes no restrictions on the mutual position
of the spectral sets $\spec(A_0)$ and $\spec(A_1)$ except that they
are disjoint and separated from each other by a distance $d$.
Sometimes, however, we will consider particular spectral dispositions
described in

\begin{hypothesis}
\label{HLgap}
Assume Hypothesis \ref{HLss}. Assume, in addition,
that either the
spectra of $A_0$ and $A_1$ are subordinated, that is,
\begin{equation}
\label{subord}
\sup\,\spec(A_0)<\inf\,\spec(A_1) \text{\, or \,} \inf\,\spec(A_0)>\sup\,\spec(A_1),
\end{equation}
or one of the sets $\spec(A_0)$ and $\spec(A_1)$ lies in a
finite gap of the other set, that is,
\begin{equation}
\label{gap}
\mathop{\rm
conv}\bigl(\spec(A_0)\bigr)\cap\spec(A_1)=\emptyset \text{\, or \,}
\spec(A_0)\cap\mathop{\rm conv}\bigl(\spec(A_1)\bigr)=\emptyset.
\end{equation}
\end{hypothesis}

Under Hypotheses \ref{HLss} or \ref{HLgap} the bound on the norm of
the inverse of the Sylvester operator \eqref{SylOp} may be given in
terms of the distance $d$ between $\spec(A_0)$ and $\spec(A_1)$. The
following result is well known.

\begin{theorem}
\label{Th-Syl-sa}
Assume Hypothesis \ref{HLss}. Let the Sylvester operator $\sS$ be
defined by \eqref{SylOp} and~\eqref{SylDom}.
\begin{enumerate}
\item[(i)] Then the inverse of $\sS$ exists and is bounded. Moreover,
the following estimate holds:
\begin{equation}
\label{Spi2}
\|\sS^{-1}\|\leq \frac{\pi}{2d}.
\end{equation}

\item[(ii)] Assume Hypothesis \ref{HLgap}. Then the following stronger
inequality holds:
\begin{equation}
\label{Sd}
\|\sS^{-1}\|\leq \frac{1}{d}.
\end{equation}
\end{enumerate}
\end{theorem}

\begin{remark}
In the generic case (i), where no assumptions on the mutual position
of the sets $\spec(A_0)$ and $\spec(A_1)$ are imposed, the existence
of a universal constant $c$ such that $\|\sS^{-1}\|\leq
\dfrac{c}{d}$ has been proven in \cite{BDM1983}.  The proof of the
fact that $c=\pi/2$ is best possible is due to R. McEachin
\cite{McE93}. For more details see \cite[Remark 2.8]{AMM}. As for
the particular spectral disposition \eqref{subord}, the bound
\eqref{Sd} is an immediate corollary to Theorem \ref{Kexp}. Since
any self-adjoint operator is simultaneously a normal operator,  in
the case of the spectral disposition \eqref{gap} the bound
\eqref{Sd} follows from Corollary \ref{Xnormal}. Sharpness of the
bound \eqref{Sd} in case (ii) is proven by an elementary example
where the spaces $\fH_0$ and $\fH_1$ are one-dimensional,
$\fH_0=\fH_1=\bbC$, and the entries $A_0=a_0$ and $A_1=a_1$ are real
numbers such that $|a_1-a_0|=d>0$.
\end{remark}

Under the assumption that both the entries $A_0$ and $A_1$ are
self-adjoint operators, below we present an existence result for the
operator Riccati equation \eqref{RicABCD}, which is written directly
in terms of the distance between the spectra of the entries $A_0$
and $A_1$ (and norms of the operators $B$ and $C$). The result is an
immediate corollary to Theorems \ref{ExistGen} and \ref{Th-Syl-sa}.
We only notice that the role of the quantity $\delta$ in the bounds
like \eqref{S1del} (see inequalities \eqref{EstLdpi} and
\eqref{EstL2} below) will be played by either $\frac{2}{\pi}d$ from
\eqref{Spi2} or $d$ from \eqref{Sd}.

\begin{theorem}
\label{ExistGenS} Assume Hypothesis \ref{HLss}.

\begin{enumerate}

\item[(i)] Then for any $B\in\cB(\fH_1,\fH_0)$ and $C\in\cB(\fH_0,\fH_1)$ such that
\begin{equation}
\label{BCpi}
\sqrt{\|B\|\|C\|}<\dfrac{d}{\pi}
\end{equation}
the Riccati equation \eqref{RicABCD} has a unique strong
solution $K$ in the ball $\cO_{d/(\pi\|B\|)}(\fH_0,\fH_1)$. This solution
satisfies the estimate
\begin{equation}
\label{EstLdpi}
\|K\|\leq \dfrac{\|C\|}{\frac{d}{\pi}+\sqrt{\frac{d^2}{\pi^2}-
\|B\|\,\|C\|}}.
\end{equation}

\item[(ii)] If the conditions of Hypothesis \ref{HLgap} also hold
then the Riccati equation \eqref{RicABCD} has a unique strong
solution $K$ in the ball $\cO_{d/(2\|B\|)}(\fH_0,\fH_1)$ whenever
$B\in\cB(\fH_1,\fH_0)$ and $C\in\cB(\fH_0,\fH_1)$ satisfy the bound
\begin{equation}
\label{BC12}
\sqrt{\|B\|\|C\|}<\dfrac{d}{2}.
\end{equation}
The solution $K$ satisfies the estimate
\begin{equation}
\label{EstL2}
\|K\|\leq \frac{\|C\|}{\frac{d}{2}+\sqrt{\frac{d^2}{4}-
\|B\|\,\|C\|}}.
\end{equation}
\end{enumerate}
\end{theorem}
\begin{remark}
The part (i) is a refinement of Theorem 3.6 in \cite{AMM} that
only claimed the existence of a weak (but not strong) solution to the Riccati
equation \eqref{RicABCD} within the ball
$\cO_{d/(\pi\|B\|)}(\fH_1,\fH_0)$.  The result of the part (ii) is new.
\end{remark}

\begin{remark}
\label{Rsep} Let $r$ be given by formula \eqref{rb} where
$\delta=\frac{2}{\pi}d$ in case {\rm(i)} and $\delta=d$ in
case~{\rm(ii)}. By Remarks \ref{RZZ} and \ref{Rrb} one concludes
that the spectrum of the block operator matrix $L$ consists of a two
disjoint components $\sigma'_0=\spec(Z_0)$ and
$\sigma'_1=\spec(Z_1)$ lying in the closed $r$-neighborhoods
$O_r\bigl(\spec(A_0)\bigr)$ and $O_r\bigl(\spec(A_1)\bigr)$  of the
corresponding spectral sets $\spec(A_0)$ and $\spec(A_1)$.
\end{remark}

\begin{remark}
\label{Rsharp0} Examples \ref{Ex1} and \ref{Ex2} below show that the
bound \eqref{EstL2} is sharp in the following sense. Given a number
$d>0$ and values of the norms $\|B\|$ and $\|C\|$ satisfying
\eqref{BC12} one can always present self-adjoint (and even rank one
or two) entries $A_0$, $A_1$ and bounded $B$ and $C$ such that in
case (ii) the bound \eqref{EstL2} turns into equality. Notice that
Examples \ref{Ex1} and \ref{Ex2} serve for the spectral
dispositions \eqref{subord} and \eqref{gap}, respectively.
\end{remark}

\begin{example}
\label{Ex1} Let $\fH_0=\fH_1=\bbC$. In this case the entries $A_0$,
$A_1$, $B$ and $C$ of \eqref{RicABCD} are simply the operators of
multiplication by numbers. Set $A_0=-\frac{d}{2}$,
$A_1=\frac{d}{2}$, $B=b$, and $C=-c$ where $b,c,$ and $d$ are
positive numbers such that $\sqrt{bc}<d/2$. The Riccati equation
\eqref{RicABCD} turns into a numeric quadratic equation whose
solutions $K^{(1)}$ and $K^{(2)}$ are given by
\begin{equation}
\label{K12}
K^{(1)}=\dfrac{c}{\frac{d}{2}+\sqrt{\frac{d^2}{4}-bc}},
\quad
K^{(2)}=\dfrac{c}{\frac{d}{2}-\sqrt{\frac{d^2}{4}-bc}}.
\end{equation}
The right-hand sides of the equalities in \eqref{K12} also represent
the norms of the corresponding solutions $K^{(1)}$ and $K^{(2)}$.
Obviously, only the solution $K^{(1)}$ satisfies the bound
$\|K\|<\frac{d}{2\|B\|}$. Also notice that the eigenvalues of the
associated $2\times 2$ matrix $L$ (which is given by \eqref{L}) read
$\lambda_-=-\sqrt{d^2/4-bc}$ and $\lambda_+=-\lambda_-$. One
observes, in particular, that $\lambda_-=A_0+BK^{(1)}$.
\end{example}

\begin{example}
\label{Ex2}
Let $\fH_0=\bbC$ and $\fH_1=\bbC^2$. Assume that
$$
A_0=0,\quad
A_1=\begin{pmatrix} -d & 0\\ \phantom{-} 0 & d \end{pmatrix},\quad
B=(0\quad b),
\quad\text{and}\quad
C=\begin{pmatrix}\phantom{-}0\\-c\end{pmatrix},
$$
where $b,c$, and $d$ are positive numbers such that
$\sqrt{bc}<{d}/{2}$. In this case the Riccati equation
\eqref{RicABCD} is easily solved explicitly. It has two solutions
$K^{(1)}=\begin{pmatrix} k^{(1)}_-\\k^{(1)}_+\end{pmatrix}$ and
$K^{(2)}=\begin{pmatrix} k^{(2)}_-\\k^{(2)}_+\end{pmatrix}$ with
$k^{(1)}_-=k^{(2)}_-=0$ and
$$
k^{(1)}_+=\dfrac{c}{\frac{d}{2}+\sqrt{\frac{d^2}{4}-bc}},
\quad
k^{(2)}_+=\dfrac{c}{\frac{d}{2}-\sqrt{\frac{d^2}{4}-bc}}.
$$
Clearly, $\|B\|=b$, $\|C\|=c$, and only the solution $K^{(1)}$
belongs to the ball $\cO_{d/(2\|B\|)}(\fH_1,\fH_0)$. Its norm is
given by the equality
\begin{equation*}
\|K^{(1)}\|=\dfrac{\|C\|}{\frac{d}{2}+\sqrt{\frac{d^2}{4}-
\|B\|\,\|C\|}}.
\end{equation*}
\end{example}

\section{$J$-symmetric perturbations}
\label{SecJsym}

In this section we deal with perturbations of spectral subspaces of
a self-adjoint operator under off-diagonal $J$-self-adjoint
perturbations.

For notational setup we adopt the following hypothesis.

\begin{hypothesis}
\label{HJS} Assume that  $A_0$ and $A_1$ are self-adjoint operators
on the Hilbert spaces $\fH_0$ and $\fH_1$ with domains $\dom(A_0)$
and $\dom(A_1)$, respectively. Let $B$ be a bounded operator from
$\fH_1$ to $\fH_0$ and $C=-B^*$. Also assume that $A$ and $V$ are
operators on $\fH=\fH_0\oplus\fH_1$ given by \eqref{Adiag} and
\eqref{Voff}, respectively, and $L=A+V$ with $\dom(L)=\dom(A)$.
\end{hypothesis}

By $J$,
\begin{equation}
\label{J2}
J=\left(\begin{array}{cr} I & 0\\ 0 & -I\end{array}\right),
\end{equation}
$\bigl($cf. \eqref{Jinv}$\bigr)$ we denote a natural involution
on the Hilbert space $\fH$ associated with its orthogonal
decomposition $\fH=\fH_0\oplus\fH_1$. Subsequently introducing the
indefinite inner product
\begin{equation}
\label{IpKs}
[x,y]=(Jx,y), \quad  x,y\in{\fH},
\end{equation}
turns ${\fH}$ into a Krein space that we denote by
${\fK}$.

A (closed) subspace $\fL\subset\fK$ is called
\textit{uniformly positive} if there is $\gamma>0$ such that
\begin{equation}
\label{mups}
[x,x]\geq \gamma\;\|x\|^2 \text{\, \,for any nonzero\, } x\in\fK.
\end{equation}
The subspace $\fL$ is called \textit{maximal uniformly positive} if
it is not a subset of any other uniformly positive subspace
of $\fK$. Uniformly negative and maximal uniformly negative
subspaces of $\fK$ are defined in a similar way. The only difference
is in the replacement of \eqref{mups} by the inequality $[x,x]\leq
-\gamma\;\|x\|^2$ that should also hold for all $x\in\fK$, $x\neq 0$.
For more definitions related to the Krein spaces we refer to
\cite{AI}, \cite{LangerEnc}.

Clearly, under Hypothesis \ref{HJS} both $V$ and $L$
are $J$-self-adjoint operators on $\fH$, that is, the products $JV$
and $JL$ are self-adjoint with respect to the initial inner product
$(\cdot,\cdot)$. This means that $V$ and $L$ are
self-adjoint on the Krein space $\fK$.

The statement below provides us with a sufficient condition for a
self-adjoint block operator matrix $L$ on $\fH$ to have purely real
spectrum and to be similar to a self-adjoint operator on $\fH$.
Notice that for the particular case where the spectra of the entries
$A_0$ and $A_1$ are subordinated, say
$\sup\spec(A_0)<\inf\spec(A_1)$, closely related results may be
found in \cite[Theorem 4.1]{AdL1995} and \cite[Theorem 3.2]{MenShk}.

\begin{theorem}
\label{Lss}
Assume Hypothesis \ref{HJS}. Suppose that the Riccati equation
\begin{equation}
\label{RicABBB}
KA_0-A_1K+KBK=-B^*
\end{equation}
has a weak (and hence strong) strictly contractive solution $K$,
$\|K\|<1.$ Then:
\begin{enumerate}
\item[(i)]  The operator matrix $L$ has a purely real spectrum and it
is similar to a self-adjoint ope\-rator on $\fH$. In particular, the
following equality holds:
\begin{equation}
\label{TLam}
L=T\Lambda T^{-1},
\end{equation}
where $T$ is a bounded and boundedly invertible operator on $\fH$ given by
\begin{equation}
\label{Ws}
T=\left(\begin{array}{ll}
I & K^* \\
K & I
\end{array}\right)
\left(\begin{array}{cc}
I-K^*K & 0 \\
0 & I-KK^*
\end{array}\right)^{-1/2}
\end{equation}
and $\Lambda$ is a block diagonal self-adjoint operator on $\fH$,
\begin{equation}
\label{Lambda}
\Lambda=\diag(\Lambda_0,\Lambda_1), \quad
\dom(\Lambda)=\dom(\Lambda_0)\oplus\dom(\Lambda_1),\,
\end{equation}
whose entries
\begin{equation}
\label{L0p}
\begin{array}{ll}
\Lambda_0=(I-K^*K)^{1/2}(A_0+BK)(I-K^*K)^{-1/2}, \\
\qquad\qquad\dom(\Lambda_0)=\Ran(I-K^*K)^{1/2}\bigr|_{\dom(A_0)},
\end{array}
\end{equation}
and
\begin{equation}
\label{L01}
\begin{array}{ll}
\Lambda_1=(I-KK^*)^{1/2}(A_1-B^*K^*)(I-KK^*)^{-1/2}, \\
\qquad\qquad\dom(\Lambda_1)=\Ran(I-KK^*)^{1/2}\bigr|_{\dom(A_1)},
\end{array}
\end{equation}
are self-adjoint operators on the corresponding component Hilbert
spaces $\fH_0$ and $\fH_1$.

\item[(ii)] The graph subspaces $\fH'_0=\cG(K)$ and $\fH'_1=\cG(K^*)$
are invariant under $L$ and mutually orthogonal with respect to the
indefinite inner product \eqref{IpKs}. Moreover,
$\fK=\fH'_0[+]\fH'_1$ where the sign ``$[+]$'' stands for the
orthogonal sum in the sense of the Krein space $\fK$. The subspace
$\fH'_0$ is  maximal uniformly positive while $\fH'_1$ maximal
uniformly negative. The restrictions of $L$ onto the subspaces
$\fH'_0$ and $\fH'_1$ are $\fK$-unitary equivalent to the
self-adjoint operators $\Lambda_0$ and $\Lambda_1$, respectively.

\end{enumerate}
\end{theorem}
\begin{proof}
In the case under consideration the second Riccati equation \eqref{RicABCD1}
associated with the operator matrix $L$ reads
\begin{equation}
\label{RicABBa}
K'A_1-A_0K'-K'B^*K'=B.
\end{equation}
Thus, it simply coincides with the corresponding adjoint
\eqref{RicAdj} of the Riccati equation \eqref{RicABBB}. By Remark
\ref{RemRAdj} this means that the adjoint of $K$, $K'=K^*$, is a
weak (and hence strong) solution to \eqref{RicABBa}. Since
$\|K^*\|=\|K\|<1$, the operators $I-K^*K$ and $I-KK^*$ are strictly
positive,
\begin{equation}
\label{KKs}
I-K^*K\geq I-\|K\|^2>0\quad \text{and}\quad I-KK^*\geq I-\|K\|^2>0,
\end{equation}
and, hence, boundedly invertible. This also means that the operator
$T$ in \eqref{Ws} is well defined and bounded. In addition, by
Remark \ref{RemKK} this implies that the operator $W$ in \eqref{W}
is boundedly invertible and, consequently, the same holds for $T$.

Now notice that by Lemma \ref{L2ident} we have
\begin{equation*}
\ran\bigl({K^*K}|_{\dom(A_0)}\bigr)\subset\dom(A_0), \quad
\ran\bigl({KK^*}|_{\dom(A_1)}\bigr)\subset\dom(A_1),
\end{equation*}
and
\begin{align}
\label{Id0c}
(I-K^*K)(A_0+BK)x &=(A_0+K^*B^*)(I-K^*K)x\quad\text{for all \,}x\in\dom(A_0),\\
\label{Id1c}
(I-KK^*)(A_1-B^*K^*)y &=(A_1-KB)(I-KK^*)y\quad\text{for all \,}y\in\dom(A_1),
\end{align}
from which one easily infers that both $\Lambda_0$ and $\Lambda_1$ are
self-adjoint operators.

By using \eqref{L0p} and \eqref{L01} one expresses the operators
$Z_0=A_0+BK$, $\dom(Z_0)=\dom(A_0)$, and  $Z_1=A_1-B^*K^*$,
$\dom(Z_0)=\dom(A)$, in terms of $\Lambda_0$ and $\Lambda_1$. Then
combining the expressions obtained with equality \eqref{LWZ} from
Corollary \ref{Ldiag} we obtain formula \eqref{TLam}. The similarity
\eqref{TLam} means, in particular, that $\spec(L)$ is a Borel
subset of $\bbR$. This completes the proof of part (i).

The $J$-orthogonality of the subspaces $\fH'_0$ and $\fH'_1$ is
obvious since for any $x,y\in\fH$ of the form
\begin{equation}
\label{xk} x=x_0\oplus Kx_0,\,\, x_0\in\fH_0,\,\,\text{and}\,\,
y=K^*y_1\oplus y_1,\,\, y_1\in\fH_1,
\end{equation}
we have $ [x,y]=(Jx,y)=(x_0,K^*y_1)-(Kx_0,y_1)=0. $ Thus, the fist
two assertions of part (ii) follow from Corollary \ref{Ldiag}\;(ii).
On the other hand, \eqref{xk} yields
$\|x\|^2\leq(1+\|K\|^2)\|x_0\|^2$ and
\mbox{$\|y\|^2\leq(1+\|K\|^2)\|y_1\|^2$}, and, hence, combined
with \eqref{KKs}, it implies $[x,x]\geq \gamma\|x\|^2$ and $[y,y]\leq
-\gamma\|y\|^2$ where $\gamma={(1-\|K\|^2)}{(1+\|K\|^2)^{-1}}>0$.
This means that $\fH'_0$ and $\fH'_1$ are maximal uniformly positive
and maximal uniformly negative subspaces, respectively.

Now introduce the operators $T_0=W_0(I-K^*K)^{-1/2}$ and
$T_1=W_1(I-KK^*)^{-1/2}$ where $W_0$ and $W_1$ are given in
\eqref{W01} assuming that $K'=K^*$. Taking into account \eqref{L0p}
and \eqref{L01}, the identities \eqref{W01Z} of Corollary
\ref{Ldiag}\;(ii) then imply
\begin{equation}
\label{T01L} T_0^{-1}L|_{\fH'_0}T_0=\Lambda_0\quad\text{and}\quad
T_1^{-1}L|_{\fH'_1}T_1=\Lambda_1.
\end{equation}
Clearly, $\Ran T_0=\fH'_0$,
$\Ran T_1=\fH'_1$, $[T_0x_0,T_0y_0]=(x_0,y_0)$ for any
$x_0,y_0\in\fH_0$, and $[T_1x_1,T_1y_1]=-(x_1,y_1)$ for any
$x_1,y_1\in\fH_1$. This means that both $T_0:\,\fH_0\to\fH'_0$ and
$T_1:\,\fH_1\to\fH'_1$ are $\fK$-unitary operators. Therefore, equalities
\eqref{T01L} prove the remaining statement of part (ii).

The proof is complete.
\end{proof}

\begin{remark}
\label{Zss}
By equalities \eqref{L0p} and \eqref{L01} the self-adjoint operators
$\Lambda_0$ and $\Lambda_1$ are similar to the operators
\begin{align}
\label{Z0p}
Z_0&=A_0+BK, \,\, \dom(Z_0)=\dom(A_0), \text{\,\, and \,\,}
Z_1=A_1-B^*K^*, \,\, \dom(Z_1)=\dom(A_1)
\end{align}
respectively, and, thus,
\begin{equation}
\label{ZLZL} \spec(\Lambda_0)=\spec(Z_0) \text{\, and \,}
\spec(\Lambda_1)=\spec(Z_1).
\end{equation}
Notice that identities \eqref{Id0c} and \eqref{Id1c} imply that the
operators $Z_0$ and $Z_1$ are self-adjoint on the corresponding
Hilbert spaces $\fH_0$ and $\fH_1$ equipped with the new inner
products $ \lal
f_0,g_0\ral_{\fH_0}=\bigl((I-K^*K)f_0,g_0\bigr)_{\fH_0} \,\,\text{
and }\,\, \lal
f_1,g_1\ral_{\fH_1}=\bigl((I-KK^*)f_1,g_1\bigr)_{\fH_1}, $
respectively.
\end{remark}

\begin{remark}
\label{RKcontr} The requirement $\|K\|<1$ is sharp in the following
sense: If there is no strictly contractive solution to the Riccati
equation \eqref{RicABBB} then the operator matrix $L$ may not be
similar to a self-adjoint operator at all. This is clearly seen from
the simple example below.
\end{remark}

\begin{example}
\label{ExNS} Let $\fH_0=\fH_1=\bbC$. Set $A_0=-\frac{d}{2}$,
$A_1=\frac{d}{2}$, and $B=b$ where $b$ and $d$ are positive numbers
such that $b\geq\frac{d}{2}$. If $b>d/2$, the Riccati equation
\eqref{RicABBB} has two solutions
$X^{(1)}=\frac{d}{2b}+i\sqrt{\frac{d^2}{4b^2}-1}$ and
$X^{(2)}=\frac{d}{2b}-i\sqrt{\frac{d^2}{4b^2}-1}$. Both $X^{(1)}$
and $X^{(2)}$ are not strictly contractive since
$\|X^{(1)}\|=\|X^{(2)}\|=1$. At the same time the spectrum of the
matrix $L$ consists of the two complex eigenvalues
$\lambda_1=i\sqrt{b^2-\frac{d^2}{4}}$ and
$\lambda_2=-i\sqrt{b^2-\frac{d^2}{4}}$. If $b=\frac{d}{2}$, the
equation \eqref{RicABBB} has the only solution $X=1$. In this case
the spectrum of the matrix $L$ is real (it consists of the only
point zero) but one easily verifies by inspection that the only
eigenvalue of $L$ has a nontrivial Jordan chain and, thus, $L$ is
not diagonalizable. Therefore, in both cases $b>d/2$ and $b=d/2$ the
matrix $L$ cannot be made similar to a self-adjoint operator.
\end{example}

The next assertion represents a quite elementary corollary to
Theorem \ref{Lss}.

\begin{lemma}
\label{Tcomp} Let the assumptions of Theorem \ref{Lss} hold. Assume,
in addition, that the spectra $\sigma'_0=\spec(Z_0)$ and
$\sigma'_1=\spec(Z_1)$ of the operators $Z_0$ and $Z_1$ given by
\eqref{Z0p} are disjoint, that is,
$\sigma'_0\cap\sigma'_1=\emptyset$. Then $\sigma'_0$ and $\sigma'_1$
are complementary spectral subsets of the block operator matrix $L$,
$\spec(L)=\sigma'_0\cup\sigma'_1$, and the graphs $\fH'_0=\cG(K)$
and $\fH'_1=\cG(K^*)$ are the spectral subspaces associated with the
subsets $\sigma'_0$ and $\sigma'_1$, respectively.
\end{lemma}
\begin{proof} By the assumption the spectra
$\spec(\Lambda_0)=\spec(Z_0)=\sigma'_0$ and
$\spec(\Lambda_1)=\spec(Z_1)=\sigma'_1$ (see Remark \ref{Zss}) of
the self-adjoint operators $\Lambda_0$ and $\Lambda_1$ given by
\eqref{L0p}, \eqref{L01} are disjoint. Hence, the spectral
projections $\sE_\Lambda(\sigma'_0)$ and $\sE_\Lambda(\sigma'_1)$ of
the self-adjoint diagonal block operator matrix
$\Lambda=\diag(\Lambda_0,\Lambda_1)$ associated with its spectral
subsets $\sigma'_0$ and $\sigma'_1$ read simply as
$$
\sE_\Lambda(\sigma'_0)=\left(\begin{array}{cc}
I & 0 \\
0 & 0
\end{array}\right)\quad\text{and}\quad\sE_\Lambda(\sigma'_1)=\left(\begin{array}{cc}
0 & 0 \\
0 & I
\end{array}\right)
$$
By Theorem \ref{Lss}\;(i) the operator $L$ is similar to the operator
$\Lambda$. This means that the similarity transforms
$\sE_L(\sigma'_0)=T\sE_\Lambda(\sigma'_0)T^{-1}$ and
$\sE_L(\sigma'_1)=T\sE_\Lambda(\sigma'_1)T^{-1}$ of the spectral
projections $\sE_\Lambda(\sigma'_0)$ and $\sE_\Lambda(\sigma'_1)$ with $T$ given
by \eqref{Ws} represent the corresponding spectral projections of
$L$. One verifies by inspections that $\sE_L(\sigma'_0)=Q_{\cG(K)}$
and $\sE_L(\sigma'_1)=Q_{\cG(K^*)}$ where $Q_{\cG(K)}$ and
$Q_{\cG(K^*)}$ are given by \eqref{PGK} assuming that $K'=K^*$. That
is, $\sE_L(\sigma'_0)$ and $\sE_L(\sigma'_1)$ are the (oblique)
projections onto the graph subspaces $\cG(K)$ and $\cG(K^*)$,
respectively, which completes the proof.
\end{proof}
\begin{remark}
The spectral projections $\sE_L(\sigma'_0)=Q_{\cG(K)}$ and
$\sE_L(\sigma'_1)=Q_{\cG(K^*)}$ are orthogonal projections with
respect to the Krein inner product \eqref{IpKs}.
\end{remark}

{}From now on we will assume that the spectra of the entries $A_0$ and
$A_1$ are disjoint and, thus, the sets $\sigma_0=\spec(A_0)$ and
$\sigma_1=\spec(A_1)$ appear to be complementary disjoint spectral
subsets of the total self-adjoint operator $A$. In such a case for
any bounded perturbation $V$ satisfying the bound $\|V\|<d/2$,
$d=\dist(\sigma_0,\sigma_1)$, the spectrum of the perturbed operator
$L=A+V$ consists of two disjoint subsets $\sigma'_0$ and
$\sigma'_1$, lying in the closed $\|V\|$-neighborhoods
$O_{\|V\|}(\sigma_0)$ and $O_{\|V\|}(\sigma_1)$ of the spectral sets
$\sigma_0=\spec(A_0)$ and $\sigma_1=\spec(A_1)$, respectively. One
can think of the sets $\sigma'_0$ and $\sigma'_1$ as the result of
the perturbation of the corresponding spectral sets $\sigma_0$ and
$\sigma_1$.

Provided that the perturbation $V$ is $J$-symmetric and $\|V\|<d/2$,
Theorem \ref{Tpi} below gives sufficient \textit{a priori}
conditions for the perturbed operator $L=A+V$ to remain similar to a
self-adjoint operator. Hence, this theorem also gives sufficient
conditions for the perturbed spectral sets $\sigma'_0$ and
$\sigma'_1$ to remain on the real axis. Furthermore, the theorem
presents the main result of the section giving for such $V$ an
\textit{a priori} norm bound on variation of the spectral subspaces
of $A$ associated with the disjoint spectral subsets $\sigma_0$ and
$\sigma_1$.

\begin{theorem}
\label{Tpi} Assume Hypothesis \ref{HJS} and choose one of the
following:
\begin{enumerate}
\item[(i)] Assume \eqref{AAdist} and set $\delta=\frac{2}{\pi}d$;
\item[(ii)] Assume \eqref{subord} or \eqref{gap} and set $\delta=d$.
\end{enumerate}
Also suppose that
\begin{equation}
\label{Bdel}
\|V\|<\frac{\delta}{2}.
\end{equation}
Then the spectrum of the operator $L$ is purely real and consists
of two disjoint components $\sigma'_0$ and $\sigma'_1$ such that
\begin{equation}
\label{sigp}
\sigma'_0\subset O_r\bigl(\spec(A_0)\bigr)\text{ \, and \, }
\sigma'_1\subset O_r\bigl(\spec(A_1)\bigr),
\end{equation}
where
$$
r=\|V\|\tanh\left(\frac{1}{2}\mathop{\rm arctanh}\frac{2\|V\|}{\delta}\right)<\|V\|.
$$
Moreover, the operator $L$ is similar to a self-adjoint operator and
the same is true for the parts of $L$ associated with the spectral
subsets $\sigma'_0$ and $\sigma'_1$. Furthermore, the following
bound holds:
\begin{equation}
\label{Thbpi} \tan\Theta_0\leq\tanh\left(\frac{1}{2}\mathop{\rm
arctanh}\frac{2\|V\|}{\delta}\right),
\end{equation}
where $\Theta_0=\Theta(\fH_0,\fH'_0)$ denotes the operator angle between the
subspace $\fH_0$ and the spectral subspace $\fH'_0$ of $L$ associated with
the spectral subset $\sigma'_0$. Exactly the
same bound holds for the operator angle $\Theta_1=\Theta(\fH_1,\fH'_1)$ between
the subspace $\fH_1$ and the spectral subspace $\fH'_1$ of $L$ associated
with the spectral subset $\sigma'_1$.
\end{theorem}
\begin{proof}
Under either assumption (i) or (ii) from Theorem \ref{ExistGenS} it
follows that the Riccati equation \eqref{RicABBB} associated with
the block operator matrix $L$ has a solution $K\in\cB(\fH_0,\fH_1)$
that is unique in the ball $\cO_{\delta/{2\|B\|}}(\fH_0,\fH_1)$ and
satisfies the bound $\bigl($see formulas \eqref{EstLdpi} and
\eqref{EstL2}$\bigr)$
\begin{equation}
\label{KbFin}
\|K\|\leq \frac{\|V\|}{\frac{\delta}{2}+\sqrt{\frac{d^2}{4}-
\|V\|^2}}
=\tanh\left(\frac{1}{2}\mathop{\rm
arctanh}\frac{2\|V\|}{\delta}\right).
\end{equation}
Here we have taken into account that $\|B\|=\|V\|$. We refer to
Remark \ref{Bhyp} regarding the use of the hyperbolic tangent in \eqref{KbFin}.

Clearly, the bound \eqref{KbFin} yields that the solution $K$ is a
strict contraction, $\|K\|<1$. Then by Theorem \ref{Lss} the block
operator matrix $L$ is similar to the self-adjoint operator
$\Lambda$ given by \eqref{Lambda}--\eqref{L01}. Hence
$\spec(L)\subset\bbR$ and $\spec(L)=\sigma'_0\cup\sigma'_1$ where
$\sigma'_0=\spec(\Lambda_0)$ and
$\sigma'_1=\spec(\Lambda_1)$. By Remark \ref{Zss} we
also have $\sigma'_0=\spec(Z_0)$ and $\sigma'_1=\spec(Z_1)$ where
$Z_0$ and $Z_1$ are given by \eqref{Z0p}. Since
$\|BK\|\leq\|V\|\|K\|\leq r$ and $\|B^*K^*\|\leq\|V\|\|K\|\leq r$,
for the spectral sets $\sigma'_0=\spec(Z_0)$ and
$\sigma'_1=\spec(Z_1)$ the inclusions \eqref{sigp} hold and these
sets are disjoint,
$\dist(\sigma'_0,\sigma'_1)\geq\delta-2r>\delta-2\|V\|>0$. To prove
the remaining statements of the theorem one only needs to apply
Lemma \ref{Tcomp} and then to notice that due to \eqref{KThet} we
have $\|\tan\Theta_0\|=\|K\|$ and hence $\tan\Theta_0\leq \|K\|$.
Similarly, \mbox{$\tan\Theta_1\leq\|K^*\|=\|K\|$}.

The proof is complete.
\end{proof}

\begin{remark}
\label{dpisharp}
By the upper continuity of the spectrum, the inclusion $\spec(L)\subset\bbR$ also
holds for $\|V\|=d/\pi$ in case (i) and for $\|V\|=d/2$ in case (ii).
\end{remark}

\begin{remark}
\label{Rsharp} In case (ii) the bounds \eqref{sigp} on the location
of $\spec(L)$ and the bound \eqref{Thbpi} on the angle $\Theta_0$
are optimal. The optimality of both \eqref{sigp} and \eqref{Thbpi}
is seen from Examples \ref{Ex1} and \ref{Ex2} where one sets $c=b$.
\end{remark}

\begin{remark}
\label{RHHKr} Under condition \eqref{Bdel} in both cases (i) and
(ii) the perturbed spectral subspaces $\fH'_0$ and $\fH'_1$ are
mutually orthogonal with respect to the Krein space inner product
\eqref{IpKs} and, thus, $\fK=\fH'_0[+]\fH'_1$. These subspaces are
maximal uniformly positive and maximal uniformly negative,
respectively. The restrictions of $L$ onto $\fH'_0$ and $\fH'_1$ are
$\fK$-unitary equivalent to the self-adjoint operators $\Lambda_0$
and $\Lambda_1$ given by \eqref{L0p} and \eqref{L01}, respectively.
By Theorem \ref{Lss}\;(ii) all this follows from the fact that
$\|K\|<1$ which we established in the proof of Theorem \ref{Tpi}.
\end{remark}

Theorem \ref{Tpi} claims that the spectrum of the block operator
matrix $L$ is purely real whenever the off-diagonal $J$-self-adjoint
perturbation $V$ satisfies the bounds $\|V\|<d/2$ in case (i) or
$\|V\|<d/\pi$ in case (ii). Recall that case (ii) corresponds to the
general spectral situation where no constraints are imposed on the
mutual positions of the spectra $\spec(A_0)$ and $\spec(A_1)$ except
for the condition \eqref{AAdist}. Now we want to prove that, in
fact, under the only condition \eqref{AAdist} the spectrum of the
operator $L$ remains purely real even if $d/\pi\leq\|V\|<d/2$, at
least in the case where the entries $A_0$ and $A_1$ are bounded. Our
proof will be based on results from \cite{LT2004-PT} and
\cite{TretterWag2003}.

\begin{theorem}
\label{Tspreal} Assume Hypothesis \ref{HJS}. Assume, in addition,
that both the entries $A_0$ and $A_1$ are bounded and such that
$\dist\bigl(\spec(A_0),\spec(A_1)\bigr)=d>0$. Also suppose that
\mbox{$\|V\|<d/2$}. Then the spectrum of the block operator matrix
$L$ is real, that is, $\spec(L)\subset\bbR$.
\end{theorem}
\begin{proof}
Under Hypothesis \ref{HLgap} and condition $\|V\|<d/2$ the inclusion
$\spec(L)\subset\bbR$ has been already proven in Theorem \ref{Tpi}
(ii). Thus, let us only consider the case that is not covered by
Hypothesis \ref{HLgap}. In this case, because of the separation
condition $\dist\bigl(\spec(A_0),\spec(A_1)\bigr)=d$, the spectrum
of $A_0$ consists of several (at least two) nonempty subsets
isolated from each other at least by the distance $2d$. Denote these
isolated spectral subsets of $A_0$ by $\sigma^{(i)}_0$,
$i=1,2,\ldots,n_0$, $n_0\geq 2$, assuming that they are numbered
from left to right (i.e. $\sup\sigma^{(i)}_0<\inf\sigma^{(i+1)}_0$),
the gap between $\sup\sigma^{(i)}_0$ and $\inf\sigma^{(i+1)}_0$
contains a nonempty subset of the spectrum of $A_1$, and
$\bigcup_{i=1}^{n_0}\sigma^{(i)}_0=\spec(A_0)$. In exactly the same
way, divide the spectrum of $A_1$ into the subsets $\sigma^{(j)}_1$,
$j=1,2,\ldots,n_1$, $n_1\geq 2$, so that
$\bigcup_{j=1}^{n_1}\sigma^{(j)}_1=\spec(A_1)$,
$\sup\sigma^{(j)}_1<\inf\sigma^{(j+1)}_1$, and
$(\sup\sigma^{(j)}_1,\inf\sigma^{(j+1)}_1)\cap\spec(A_0)\neq\emptyset$.
Denote by $\fH_0^{(i)}$, $i=1,2,\ldots,n_0$, and $\fH_1^{(j)}$,
$j=1,2,\ldots,n_1$, the spectral subspaces of the operators $A_0$
and $A_1$ associated with the corresponding spectral subsets
$\sigma_0^{(i)}$ and $\sigma_1^{(j)}$. Surely,
$\oplus_{i=0}^{n_0}\fH_0^{(i)}=\fH_0$ and
$\oplus_{i=1}^{n_1}\fH_1^{(i)}=\fH_1$.

Now take arbitrary unit vectors
\begin{equation}
\label{e01}
e^{(i)}_0\in\fH^{(i)}_0, \,\,
\|e^{(i)}_0\|=1, \,\,i=1,2,\ldots,n_0,\text{\, and \,}
e^{(j)}_1\in\fH^{(j)}_1,\,\,
\|e^{(j)}_1\|=1,\,\, j=1,2,\ldots,n_1,
\end{equation}
and construct numerical matrices $\sA_0$, $\sA_1$, and $\sB$ with the entries
\begin{equation*}
\sA_{0,ik}=( A_0 e_k^{(0)},e_i^{(0)}),\quad
\sA_{1,jl}=( A_1 e_l^{(1)},e_j^{(1)}),\quad\text{and}\quad
\sB_{ij}=( B e_j^{(1)},e_i^{(0)}),
\end{equation*}
respectively. Consider the matrices $\sA_0$ and $\sA_1$ as operators
resp. on $\widehat{\fH}_0=\bbC^{n_0}$ and
$\widehat{\fH}_1=\bbC^{n_1}$, and $\sB$ as an operator from
$\widehat{\fH}_1$ to $\widehat{\fH}_0$. Out of the matrices $\sA_0$
and $\sA_1$ construct the block diagonal matrix
$\sA=\diag(\sA_0,\sA_1)$ and out of $\sB$ and $\sB^*$ the
off-diagonal matrix $\sV=\left(\begin{array}{cr} 0 & \sB \\-\sB^* &
0\end{array}\right)$. Both matrices $\sA$ and $\sV$ have dimension
$n\times n$ where $n=n_0+n_1$, and we consider them as operators on
the $n$-dimensional space
$\widehat{\fH}=\widehat{\fH}_0\oplus\widehat{\fH}_1$.

Our nearest goal is to prove that the spectrum of the operator
$\sL=\sA+\sV$ is real. To this end, first, introduce the
indefinite inner product
\begin{equation}
\label{IpK}
[x,y]=( x_0,y_0)_{\widehat{\fH}_0} -(
x_1,y_1)_{\widehat{\fH}_1}, \quad x=x_0\oplus x_1,\,y=y_0\oplus
y_1, \,\, x_0,y_0\in\widehat{\fH}_0,\,x_1,y_1\in\widehat{\fH}_1,
\end{equation}
which turns the Hilbert space $\widehat{\fH}$ into a Krein
(Pontrjagin) space. We denote the latter by $\widehat{\fK}$. The
operator $\sA$ is self-adjoint both on $\widehat{\fH}$ and
$\widehat{\fK}$ while $\sB$ only on $\widehat{\fK}$.

Then notice that for different $i$ and $k$ the vectors $e_i^{(0)}$
and $e_k^{(0)}$ belong to the different (and mutually orthogonal)
spectral subspaces of $A_0$ and, hence,
$\sA_{0,ik}=\lambda^{(0)}_i\delta_{ik}$ where $\lambda^{(0)}_i=(
A_0 e_i^{(0)},e_i^{(0)})$ and $\delta_{ik}$ is the Kronecker's
delta. Similarly, $\sA_{1,jl}=\lambda^{(1)}_j\delta_{jl}$
where $\lambda^{(1)}_j=( A_1 e_j^{(0)},e_j^{(0)})$. Clearly,
both $\lambda_i^{(0)}$, $i=1,2,\ldots,n_0$, and $\lambda_j^{(1)}$,
$j=1,2,\ldots,n_1$, are simple eigenvalues of $\sA$ and, by
construction of $\sA$, one has
$\lambda^{(0)}_i\in\conv(\sigma_0^{(i)})$ and
$\lambda^{(1)}_j\in\conv(\sigma_1^{(j)})$. This yields
\begin{equation}
\min_{i,k,\,i\neq k}|\lambda^{(0)}_i-\lambda^{(0)}_k|\geq 2d,\quad
\min_{j,l,\,j\neq l}|\lambda^{(1)}_j-\lambda^{(1)}_l|\geq 2d,\quad\text{and}\quad
 \min_{i,j} |\lambda^{(0)}_i-\lambda^{(1)}_j|\geq d.
\end{equation}
It is also obvious that, with respect to the inner product
\eqref{IpK}, the eigenvalues $\lambda^{(0)}_i$, $i=1,2,\ldots,n_0$,
are of positive type, while the eigenvalues $\lambda^{(1)}_j$,
$j=1,2,\ldots,n_1$, are of negative type.

Now to prove the inclusion $\spec(\sL)\subset\bbR$ it only remains
to observe that $\|\sV\|\leq\|V\|<d/2$ and then to apply
\cite[Corollary 3.4]{LT2004-PT} (cf. \cite[Theorem 1.2]{Caliceti2}).

Since the inclusion $\spec(\sL)\subset\bbR$ holds for any choice of
the vectors \eqref{e01}, one then concludes that also
$W^n(L)\subset\bbR$ where $W^n(L)$ denotes the block numerical range
(see \cite[Definition 2.1]{TretterWag2003}) of the operator $L$ with
respect to the decomposition
\begin{equation}
\label{fHn}
\fH=\fH_0^{(1)}\oplus\ldots\oplus\fH_0^{(n_0)}
\oplus\fH_1^{(1)}\oplus\ldots\oplus\fH_1^{(n_1)}.
\end{equation}
By \cite[Theorem 2.5]{TretterWag2003} we have
$\spec(L)\subset\overline{W^n(L)}$. Hence, $\spec(L)\subset\bbR$,
which completes the proof.
\end{proof}
\begin{remark}
\label{d2sharp} By the upper continuity of the spectrum, under the
hypothesis of Theorem \ref{Tspreal} the spectrum of $L=A+V$ is real
also for $\|V\|=d/2$ (cf. Remark \ref{dpisharp}).
\end{remark}
\begin{remark}
\label{RemRd2} Under the assumptions of Theorems \ref{Tpi} (ii) or
\ref{Tspreal} the requirement $\|V\|\leq d/2$ guaranteeing the
inclusion $\spec(L)\subset\bbR$ is sharp. This is seen from Example
\ref{ExNS} with $b>d/2$.
\end{remark}

\section{Quantum harmonic oscillator under a $\cP\cT$-symmetric perturbation}
\label{SecExHO}

Let $A$ be the Schr\"odinger operator for a one-dimensional quantum
harmonic oscillator (see, e.g., \cite[Chapter 12]{MessahI}). The
corresponding Hilbert space is $\fH=L_2(\bbR)$. Assuming that the
units are chosen in such a way that $\hbar=m=\omega=1$, the operator
$A$ reads
\begin{align}
\label{Aho}
(Af)(x)&=-\frac{1}{2}\frac{d^2}{dx^2}f(x)+\frac{1}{2}x^2 f(x),\quad
\dom(A)=\biggl\{f\in W^2_2(\bbR)\,\,\biggl|\,\,
\int_\bbR dx\; x^4|f(x)|^2<\infty\biggr\},\,\,
\end{align}
where $W_2^2(\bbR)$ denotes the Sobolev space of those
$L_2(\bbR)$-functions that have their second derivatives in
$L_2(\bbR)$. The subspaces
\begin{equation}
\label{Hho}
\fH_0=L_{2,\textrm{even}}(\bbR) \text{\, and \,}
\fH_1=L_{2,\textrm{odd}}(\bbR)
\end{equation}
of even and odd functions are the
spectral subspaces of the (self-adjoint) operator $A$ associated
with the spectral subsets
$$
\sigma_0=\spec(A\bigl|_{\fH_0})=\{n+1/2\,\,\bigl|\,\, n=0,2,4,\dots\}
\text{\, and \,}\sigma_1=\spec(A\bigl|_{\fH_1})=\{n+1/2\,\,\bigl|\,\,n=1,3,5\ldots\},
$$
respectively (see, e.g., \cite[p.\,142]{RS-I}). Clearly,
$\fH=\fH_0\oplus\fH_1$, the spectral sets $\sigma_0$ and $\sigma_1$
are dis\-joint,
\begin{equation}
\label{d1}
d=\dist(\sigma_0,\sigma_1)=1, \text{\, and
\,}\sigma_0\cup\sigma_1=\spec(A).
\end{equation}

Let $\cP$ be the parity operator on $L_2(\bbR)$, $(\cP
f)(-x)=f(-x)$, and $\cT$ the (antilinear) operator of complex
conjugation, $(\cT f)(x)=\overline{f(x)}$, $f\in L_2(\bbR)$. An
operator $V$ on $L_2(\bbR)$ is called $\cP\cT$-symmetric if it
commutes with the product $\cP\cT$, that is, $\cP\cT V =V\cP\cT$
(see, e.g. \cite{Caliceti1,Caliceti2} and references therein).

In a particular case where the $\cP\cT$-symmetric potential $V$ is
an operator of mul\-ti\-pli\-cation by a function $V(\cdot)$ of
$L_\infty(\bbR)$, the following equality holds (see, e.g.,
\cite{AlbFK}; cf. \cite{LT2004-PT}):
\begin{equation}
\label{VPT}
\overline{V(x)}=V(-x) \text{\, for a.e. }x\in\bbR
\end{equation}
and hence
\begin{equation}
\label{PVP}
V^*=\cP V \cP.
\end{equation}
Observe that the parity operator $\cP$ represents nothing but the
involution \eqref{Jinv} associated with the complementary spectral
subspaces \eqref{Hho} of the oscillator Hamiltonian \eqref{Aho}.
Therefore, the equality \eqref{PVP} implies that the
$\cP\cT$-symmetric multiplication operator $V$ is $J$-self-adjoint
with respect the involution $J=\cP$.

Any bounded complex-valued function $V$ on $\bbR$ possessing the
property \eqref{VPT} admits the representation
\begin{equation}
\label{VPT1}
V(x)=a(x)+i b(x)
\end{equation}
where both $a$ and $b$ are real-valued functions
such that
$$
a(-x)=a(x)\text{\, and \,}b(-x)=-b(x)
\text{\, for any \,}x\in\bbR.
$$
The terms $V_\mathrm{diag}(x)=a(x)$ and $V_\mathrm{off}(x)=ib(x)$
represent the corresponding parts of the multiplication operator $V$
that are diagonal and off-diagonal with respect to the orthogonal
decomposition $\fH=\fH_0\oplus\fH_1$, that is, with respect to the
decomposition $L_2(\bbR)=L_{2,\textrm{even}}(\bbR)\oplus
L_{2,\textrm{odd}}(\bbR)$.

Now assume that $V$ is an arbitrary bounded off-diagonal operator on
$\fH=L_2(\bbR)$ being $J$-self-adjoint with respect to the
involution $J=\cP$. One can choose in particular a $\cP\cT$-symmetric
potential \eqref{VPT1} with $a=0$. By taking into account
\eqref{d1}, from \cite[Theorem 1.2]{Caliceti1} it follows that the
spectrum of the perturbed oscillator Hamiltonian $L=A+V$,
$\dom(L)=\dom(A)$, remains real (and discrete) whenever $\|V\|\leq
1/2$. If, in addition, the bound $\|V\|<1/\pi$ is satisfied then one
can tell much more: Under such a bound Theorem \ref{Tpi}\;(i) implies
that $L$ is similar to a self-adjoint operator. This theorem also
gives bounds on the variation of the spectral subspaces \eqref{Hho}:
\begin{equation*}
\tan\Theta_j\leq\tanh\left(\frac{1}{2}\mathop{\rm
arctanh}(\pi\|V\|)\right)<1,\quad j=0,1,
\end{equation*}
where $\Theta_j=\Theta(\fH_j,\fH'_j)$ stands for the operator angle
between the subspace $\fH_j$ and the spectral subspace $\fH'_j$ of
the perturbed oscillator Hamiltonian $L=A+V$ associated with the
spectral subset $\sigma'_j=\spec(L)\cap O_{\|V\|}(\sigma_j)$,
$j=0,1.$

\vspace*{2mm} \noindent {\bf Acknowledgments.} The authors
thank S.\,M. Fei for his useful remarks on $\cP\cT$-sym\-met\-ric
operators. A. K. Motovilov and A. A. Shkalikov gratefully acknowledge
the kind hospitality of the Institut f\"ur Angewandte Mathematik,
Universit\"at Bonn, where the main part of this research has been
performed.


\end{document}